\documentclass[11pt,a4paper]{article}
\usepackage[T1]{fontenc}
\usepackage{lmodern,microtype}
\usepackage{amsmath,amssymb,amsthm,mathtools}
\usepackage[margin=2.65cm,headheight=15pt]{geometry}
\usepackage{xcolor,enumitem,booktabs,fancyhdr}
\usepackage[colorlinks=true,linkcolor=blue!45!black,citecolor=blue!45!black,urlcolor=blue!45!black]{hyperref}
\hypersetup{pdftitle={Metric Correspondences and Low-Rank Deformations for Conformal and m-Modified Conformal Vector Fields}}
\newtheorem{theorem}{Theorem}[section]
\newtheorem{lemma}[theorem]{Lemma}
\newtheorem{corollary}[theorem]{Corollary}

\theoremstyle{definition}

\newtheorem{example}[theorem]{Example}
\theoremstyle{remark}
\newtheorem{remark}[theorem]{Remark}
\numberwithin{equation}{section}
\DeclareMathOperator{\rank}{rank}

\DeclareMathOperator{\Scal}{Scal}
\DeclareMathOperator{\Ric}{Ric}
\setlist{nosep,leftmargin=2em}
\allowdisplaybreaks[2]
\title{Metric Correspondences and Low-Rank Deformations\\ for Conformal and $m$-Modified Conformal Vector Fields}
\author{}
\date{}
\begin{document}
\maketitle
\begin{abstract}
We fix a smooth vector field and a positive parameter $m$, and compare the nondegenerate metrics for which the field is conformal or $m$-modified conformal with the same prescribed smooth factor. Metrics of arbitrary signature are allowed. Passing to inverse metrics gives a linear transport equation and an explicit correspondence on the exact nondegeneracy domains. At a regular point, transverse initial data yield a signature-preserving correspondence of metric germs. For a nowhere-zero field and a Killing background in dimension at least three, every metric satisfying the modified conformal equation and differing from the background by a tensor of rank at most one has factor zero; its inverse-metric difference is $rX\otimes X$ for a unique smooth function $r$. Under a complete flow, such global deformations exist precisely when the field is null for the background and admits a smooth function $\tau$ with $X\tau=1$. An auxiliary field satisfying $[X,Z]=X$ gives a second correspondence under orthogonality assumptions. Finally, for a nowhere-zero complete field, a Killing metric and a zero-factor modified conformal metric have inverse-metric difference of rank at least two whenever that rank is constant. In the rank-two case, the pair determines a unique auxiliary field. Examples distinguish the roles of dimension, completeness, constant rank, the zero-factor condition, and orthogonality.
\end{abstract}


\section{Introduction}\label{sec:introduction}

Conformal vector fields describe infinitesimal changes of a metric by pointwise scaling. A quadratic correction involving the metric dual of the field leads to the $m$-modified conformal equation
\begin{equation}\label{eq:modified-intro}
\frac12\mathcal L_Xg-\frac1m g(X,\cdot)\otimes g(X,\cdot)=cg.
\end{equation}
Here $m>0$ is fixed and $c$ is a smooth function. We study this equation with $X$ and $c$ prescribed and the metric $g$ unknown. In particular, we ask how its solutions are related to metrics $B$ satisfying the ordinary conformal equation
\begin{equation}\label{eq:conformal-intro}
\frac12\mathcal L_XB=cB
\end{equation}
for the same field and the same factor. The distinction between fixing a field and varying a metric is essential: the metric dual in \eqref{eq:modified-intro} changes with the unknown metric.

Zhang and Chen introduced $m$-modified conformal vector fields and studied constant-length fields and left-invariant fields on Riemannian and pseudo-Riemannian Lie groups \cite{ZC}. Subsequent work has developed triviality results in the Riemannian setting. Poddar and Sharma proved, in particular, that a compact Riemannian manifold admits no nontrivial $m$-modified homothetic vector field \cite[Theorem~2.1]{PS}. Deshmukh, Bin Turki, Chen, and Fatima obtained further triviality results for general and gradient fields under additional compactness, completeness, or curvature hypotheses \cite{DBTC}. In a different direction, Li, Deng, and Zhang gave a criterion for the existence of a nonzero left-invariant field for some nondegenerate metric on an even-dimensional Lie group, and used a related algebraic metric shear \cite[Theorem~A and Lemma~4.6]{LDZ}. These results motivate a comparison that keeps the field fixed and describes the admissible metrics and their differences directly.

Equation \eqref{eq:modified-intro} also belongs to a broader family of geometric equations. Nurowski and Randall considered
\[
\frac12\mathcal L_Xg+c_1g(X,\cdot)\otimes g(X,\cdot)-c_2\Ric_g=\lambda g
\]
with constant parameters \cite[Section~1]{NR}. For constant $c$, taking $c_1=-1/m$, $c_2=0$, and $\lambda=c$ gives \eqref{eq:modified-intro}. Their constant-parameter formulation does not by itself cover an arbitrary prescribed function $c$. The quasi-Yamabe equation of Chen and Deshmukh is
\[
\frac12\mathcal L_Xg=(\Scal_g-\lambda)g+
\mu\,g(X,\cdot)\otimes g(X,\cdot)
\]
\cite[Equation~(1.3)]{CD}. Its specialization $\mu=1/m$ has the form \eqref{eq:modified-intro}, with $c=\Scal_g-\lambda$. For a prescribed $c$, a solution of \eqref{eq:modified-intro} is a quasi-Yamabe soliton of this type only if $\Scal_g-c$ is constant. No such scalar-curvature condition is imposed here.

Low-rank changes of indefinite metrics have another established source in generalized Kerr--Schild geometry. In Lorentzian signature, these changes have the form $\widetilde g=g+2H\ell\otimes\ell$, with $\ell$ null. Coll, Hildebrandt, and Senovilla studied the associated infinitesimal symmetries through equations of the form $\mathcal L_\xi g=2h\ell\otimes\ell$ and $\mathcal L_\xi\ell=a\ell$ \cite[Section~2]{CHS}. The rank-one null deformations below have this algebraic form in Lorentzian signature. Our classification starts with an unrestricted rank bound on the metric difference; its direction is a conclusion of the field equation. We allow arbitrary signature and keep the nondegeneracy and global flow conditions explicit.

The basic observation is that \eqref{eq:modified-intro} becomes linear after passing to the inverse metric:
\begin{equation}\label{eq:linear-intro}
\mathcal L_X(g^{-1})+2cg^{-1}=-\frac2m X\otimes X.
\end{equation}
The ordinary conformal equation gives the corresponding homogeneous equation. Thus a translation of inverse metrics by $rX\otimes X$ relates the two equations whenever
\[
\left(Xr+2cr+\frac2m\right)X\otimes X=0.
\]
The translation solves the tensor equation, but a metric correspondence also requires an exact test for invertibility and a signature calculation. The relevant inverse-update identity is classical; see Hager \cite[p.~221]{Hager}. We give its proof and the required signature analysis before applying it. Local solvability near a regular point does not settle the global problem: the scalar transport equation may have no global solution, and a transported inverse metric may become degenerate.

The results are organized around these two obstructions and the rank of the deformation. Section~\ref{sec:transport} gives mutually inverse transformations on explicitly specified nondegeneracy domains. It also treats nonzero constant factors and gives a signature-preserving correspondence between all metric germs at a regular point, normalized by agreement on a transverse hypersurface. The treatment of zeros of $X$ retains the tensor transport condition and does not assume that a scalar flowbox argument extends across a zero.

Section~\ref{sec:rank-one} assumes that $X$ is nowhere zero and that $B$ is a Killing metric for $X$. In dimension $n\geq3$, the condition $\rank(g-B)\leq1$ forces $c=0$ and
\[
g^{-1}=B^{-1}+rX\otimes X,\qquad
Xr=-\frac2m,\qquad 1+rB(X,X)\neq0.
\]
The three-dimensional case requires a separate argument involving the radical of a restricted bilinear form and the metric volume density. If $X$ is complete, a global deformation exists precisely when $B(X,X)=0$ and there is a smooth function $\tau$ with $X\tau=1$. All such deformations are
\[
g=B+\frac{2\tau}{m}B(X,\cdot)\otimes B(X,\cdot).
\]
The time-function condition is equivalent to a global product decomposition that carries $X$ to translation on the $\mathbb R$ factor. In dimension two these conclusions hold when the zero factor is imposed in advance; an example shows why this restriction matters.

Section~\ref{sec:auxiliary-field} starts with an auxiliary field $Z$ satisfying $[X,Z]=X$. Under $B(X,X)=B(X,Z)=0$, it constructs a signature-preserving correspondence that also preserves the full metric dual of $X$. Section~\ref{sec:constant-rank} starts instead with two already nondegenerate metrics: a Killing metric $B$ and a zero-factor modified conformal metric $g$. For a nowhere-zero complete field $X$, a constant-rank difference $g^{-1}-B^{-1}$ has rank at least two. Rank two is characterized by a unique smooth field $Z$, pointwise independent of $X$, such that
\[
g^{-1}-B^{-1}=-\frac1m(X\otimes Z+Z\otimes X),\qquad [X,Z]=X.
\]
This characterization requires neither of the orthogonality conditions used for the construction, and it need not preserve the metric dual of $X$. The examples identify the limits of each statement.

Throughout, $M$ is a connected smooth Hausdorff second-countable manifold without boundary, of dimension $n\geq2$. A metric is a smooth real nondegenerate symmetric covariant two-tensor. Its ordered signature $(p,q)$ records the numbers of positive and negative directions. We compare metrics themselves, without taking a quotient by diffeomorphisms or requiring invariance under an additional group action. Completeness always means completeness of the vector field $X$, rather than geodesic completeness of a metric. Section~\ref{sec:transport} permits zeros of $X$; Sections~\ref{sec:rank-one}--\ref{sec:constant-rank} assume that $X$ is nowhere zero.

The inverse metric $g^{-1}$ is a contravariant two-tensor. For one-forms $\eta,\xi$ and vector fields $U,V$, we use
\[
(U\otimes V)(\eta,\xi)=\eta(U)\xi(V),\qquad
(\eta\otimes\xi)(U,V)=\eta(U)\xi(V).
\]
Every symmetric product is written as the sum of its two tensor products, with no implicit factor of $1/2$. Tensor rank refers to the associated map $TM\to T^*M$ or $T^*M\to TM$, according to variance. For a metric $B$, the notation $B^\sharp\alpha$ denotes the unique vector satisfying $B(B^\sharp\alpha,\cdot)=\alpha$. All coordinate summations are written explicitly.

\section{Transport and correspondences between nondegenerate metrics}
\label{sec:transport}

Fix a smooth vector field $X$, a smooth function $c\in C^\infty(M)$, and a constant $m>0$. Zeros of $X$ are allowed throughout this section. We write
\begin{align*}
\mathcal C_c(X)
&=\left\{B:\frac12\mathcal L_XB=cB\right\},\\
\mathcal M_c(X;m)
&=\left\{g:\frac12\mathcal L_Xg
-\frac1m g(X,\cdot)\otimes g(X,\cdot)=cg\right\},
\end{align*}
where all tensors $B$ and $g$ in these sets are smooth, real, symmetric, and nondegenerate. Their signatures are unrestricted.

\begin{lemma}[Inverse metrics and contraction identities]
\label{lem:inverse-equation}
For any metric $g$,
\begin{equation}
g\in\mathcal M_c(X;m)
\quad\Longleftrightarrow\quad
\mathcal L_X(g^{-1})+2cg^{-1}=-\frac2m X\otimes X.
\label{eq:inverse-modified}
\end{equation}
Every such metric also satisfies
\begin{align}
\mathcal L_X\bigl(g(X,\cdot)\bigr)
&=2\left(c+\frac{g(X,X)}m\right)g(X,\cdot),\notag\\
X\bigl(g(X,X)\bigr)
&=2c\,g(X,X)+\frac2m g(X,X)^2,\notag\\
\operatorname{div}_gX&=nc+\frac{g(X,X)}m.
\label{eq:contraction-identities}
\end{align}
The divergence is defined using the positive metric volume density, so no orientation is required. Similarly,
\[
B\in\mathcal C_c(X)
\quad\Longleftrightarrow\quad
\mathcal L_X(B^{-1})+2cB^{-1}=0.
\]
In particular, if $\mathcal L_XB=0$, then
\[
\mathcal L_X\bigl(B(X,\cdot)\bigr)=0,\qquad
X\bigl(B(X,X)\bigr)=0,\qquad
\operatorname{div}_BX=0.
\]
\end{lemma}

\begin{proof}
In local coordinates, the inverse-metric identity is
\[
\sum_{k=1}^n g^{ik}g_{kj}=\delta^i{}_j.
\]
Taking its Lie derivative gives
\[
\sum_{k=1}^n(\mathcal L_Xg^{-1})^{ik}g_{kj}
+\sum_{k=1}^n g^{ik}(\mathcal L_Xg)_{kj}=0.
\]
Multiply by $g^{j\ell}$ and sum over $j$. The contractions are
\begin{align*}
0
&=\sum_{j=1}^n\sum_{k=1}^n
(\mathcal L_Xg^{-1})^{ik}g_{kj}g^{j\ell}
+\sum_{j=1}^n\sum_{k=1}^n
g^{ik}(\mathcal L_Xg)_{kj}g^{j\ell}\\
&=\sum_{k=1}^n(\mathcal L_Xg^{-1})^{ik}\delta_k{}^\ell
+\sum_{j=1}^n\sum_{k=1}^n
g^{ik}(\mathcal L_Xg)_{kj}g^{j\ell}\\
&=(\mathcal L_Xg^{-1})^{i\ell}
+\sum_{j=1}^n\sum_{k=1}^n
g^{ik}(\mathcal L_Xg)_{kj}g^{j\ell}.
\end{align*}
Consequently,
\[
(\mathcal L_Xg^{-1})^{ij}
=-\sum_{k=1}^n\sum_{\ell=1}^n
g^{ik}(\mathcal L_Xg)_{k\ell}g^{\ell j}.
\]
In components, the modified conformal equation reads
\[
(\mathcal L_Xg)_{k\ell}
=2cg_{k\ell}+\frac2m\sum_{a=1}^n\sum_{b=1}^n
g_{ka}X^ag_{\ell b}X^b.
\]
Substitution into the inverse-metric formula gives
\begin{align*}
(\mathcal L_Xg^{-1})^{ij}
&=-\sum_{k=1}^n\sum_{\ell=1}^n g^{ik}
\left(2cg_{k\ell}
+\frac2m\sum_{a=1}^n\sum_{b=1}^n
g_{ka}X^ag_{\ell b}X^b\right)g^{\ell j}\\
&=-2c\sum_{k=1}^n\sum_{\ell=1}^n
g^{ik}g_{k\ell}g^{\ell j}\\
&\quad-\frac2m\sum_{a=1}^n\sum_{b=1}^n
\left(\sum_{k=1}^n g^{ik}g_{ka}\right)X^a
\left(\sum_{\ell=1}^n g_{\ell b}g^{\ell j}\right)X^b\\
&=-2c\sum_{\ell=1}^n\delta^i{}_\ell g^{\ell j}
-\frac2m\sum_{a=1}^n\sum_{b=1}^n
\delta^i{}_aX^a\delta^j{}_bX^b\\
&=-2cg^{ij}-\frac2m X^iX^j.
\end{align*}
Conversely, suppose the last identity holds. Multiplication by $g$ in both slots gives
\begin{align*}
(\mathcal L_Xg)_{ij}
&=-\sum_{k=1}^n\sum_{\ell=1}^n
g_{ik}(\mathcal L_Xg^{-1})^{k\ell}g_{\ell j}\\
&=-\sum_{k=1}^n\sum_{\ell=1}^n
g_{ik}\left(-2cg^{k\ell}-\frac2mX^kX^\ell\right)g_{\ell j}\\
&=2c\sum_{k=1}^n\sum_{\ell=1}^n g_{ik}g^{k\ell}g_{\ell j}
+\frac2m\sum_{k=1}^n\sum_{\ell=1}^n g_{ik}X^kg_{\ell j}X^\ell\\
&=2c\sum_{\ell=1}^n\delta_i{}^\ell g_{\ell j}
+\frac2m\left(\sum_{k=1}^n g_{ik}X^k\right)
\left(\sum_{\ell=1}^n g_{j\ell}X^\ell\right)\\
&=2cg_{ij}+\frac2m g(X,\partial_{x^i})g(X,\partial_{x^j}).
\end{align*}
Thus $\mathcal L_Xg=2cg+(2/m)g(X,\cdot)\otimes g(X,\cdot)$, which proves \eqref{eq:inverse-modified} in both directions.

Since $[X,X]=0$, contraction commutes with $\mathcal L_X$ when the contracted vector is $X$. Thus
\begin{align*}
\mathcal L_X\bigl(g(X,\cdot)\bigr)
&=(\mathcal L_Xg)(X,\cdot)+g([X,X],\cdot)\\
&=\left(2cg+\frac2m g(X,\cdot)\otimes g(X,\cdot)\right)(X,\cdot)\\
&=2c\,g(X,\cdot)+\frac2m g(X,X)g(X,\cdot)\\
&=2\left(c+\frac{g(X,X)}m\right)g(X,\cdot).
\end{align*}
Contracting this identity once more gives
\begin{align*}
X\bigl(g(X,X)\bigr)
&=\bigl(\mathcal L_X(g(X,\cdot))\bigr)(X)+g(X,[X,X])\\
&=\left(2c\,g(X,\cdot)+\frac2m g(X,X)g(X,\cdot)\right)(X)\\
&=2c\,g(X,X)+\frac2m g(X,X)^2.
\end{align*}

Write $\rho=\sqrt{|\det(g_{ij})|}$ for the coefficient of the positive metric volume density. Differentiation of the determinant gives
\begin{align*}
\partial_k\rho
&=\frac\rho2\,\partial_k\log|\det(g_{ij})|
=\frac\rho2\sum_{i=1}^n\sum_{j=1}^n
g^{ij}\partial_kg_{ij},\\
\operatorname{div}_gX
&=\frac1\rho\sum_{k=1}^n\partial_k(\rho X^k)\\
&=\sum_{k=1}^n\partial_kX^k
+\frac12\sum_{i=1}^n\sum_{j=1}^n\sum_{k=1}^n
g^{ij}X^k\partial_kg_{ij}.
\end{align*}
The coordinate formula for the Lie derivative gives the same expression:
\begin{align*}
\frac12\sum_{i=1}^n\sum_{j=1}^n g^{ij}(\mathcal L_Xg)_{ij}
&=\frac12\sum_{i=1}^n\sum_{j=1}^n\sum_{k=1}^n g^{ij}
\left(X^k\partial_kg_{ij}+g_{kj}\partial_iX^k+g_{ik}\partial_jX^k\right)\\
&=\frac12\sum_{i=1}^n\sum_{j=1}^n\sum_{k=1}^n
g^{ij}X^k\partial_kg_{ij}\\
&\quad+\frac12\sum_{i=1}^n\sum_{k=1}^n\delta^i{}_k\partial_iX^k
+\frac12\sum_{j=1}^n\sum_{k=1}^n\delta^j{}_k\partial_jX^k\\
&=\frac12\sum_{i=1}^n\sum_{j=1}^n\sum_{k=1}^n
g^{ij}X^k\partial_kg_{ij}+\sum_{k=1}^n\partial_kX^k\\
&=\operatorname{div}_gX.
\end{align*}
Substituting the modified conformal equation into this trace gives
\begin{align*}
\operatorname{div}_gX
&=\frac12\sum_{i=1}^n\sum_{j=1}^n g^{ij}
\left(2cg_{ij}+\frac2m\sum_{a=1}^n\sum_{b=1}^n
g_{ia}X^ag_{jb}X^b\right)\\
&=c\sum_{i=1}^n\sum_{j=1}^n g^{ij}g_{ij}
+\frac1m\sum_{a=1}^n\sum_{b=1}^n\sum_{j=1}^n
\left(\sum_{i=1}^n g^{ij}g_{ia}\right)X^ag_{jb}X^b\\
&=nc+\frac1m\sum_{a=1}^n\sum_{b=1}^n\sum_{j=1}^n
\delta^j{}_aX^ag_{jb}X^b\\
&=nc+\frac1m\sum_{a=1}^n\sum_{b=1}^n g_{ab}X^aX^b
=nc+\frac{g(X,X)}m.
\end{align*}

For an ordinary conformal metric $B$, inverse differentiation yields
\begin{align*}
(\mathcal L_XB^{-1})^{ij}
&=-\sum_{k=1}^n\sum_{\ell=1}^n B^{ik}(2cB_{k\ell})B^{\ell j}\\
&=-2c\sum_{\ell=1}^n\delta^i{}_\ell B^{\ell j}
=-2cB^{ij}.
\end{align*}
Conversely, if $\mathcal L_XB^{-1}=-2cB^{-1}$, then
\begin{align*}
(\mathcal L_XB)_{ij}
&=-\sum_{k=1}^n\sum_{\ell=1}^n
B_{ik}(-2cB^{k\ell})B_{\ell j}\\
&=2c\sum_{\ell=1}^n\delta_i{}^\ell B_{\ell j}
=2cB_{ij}.
\end{align*}
The corresponding contractions, with the contraction steps displayed, are
\begin{align*}
\mathcal L_X\bigl(B(X,\cdot)\bigr)
&=(\mathcal L_XB)(X,\cdot)+B([X,X],\cdot)
=(2cB)(X,\cdot)=2cB(X,\cdot),\\
X\bigl(B(X,X)\bigr)
&=(\mathcal L_XB)(X,X)+B([X,X],X)+B(X,[X,X])\\
&=(2cB)(X,X)=2cB(X,X),\\
\operatorname{div}_BX
&=\frac12\sum_{i=1}^n\sum_{j=1}^n B^{ij}(2cB_{ij})
=c\sum_{i=1}^n\delta^i{}_i=nc.
\end{align*}
For $c=0$, these formulas become
\[
\mathcal L_X\bigl(B(X,\cdot)\bigr)=0,\qquad
X\bigl(B(X,X)\bigr)=0,\qquad \operatorname{div}_BX=0,
\]
as asserted.
\end{proof}

The inversion formula in the following lemma is the symmetric rank-one case of the classical inverse-update formula; see \cite[p.~221]{Hager}. We include its proof, together with the signature criterion needed below.

\begin{lemma}[Rank-one inversion and signature]
\label{lem:rank-one-inversion}
Let $B$ be a nondegenerate symmetric bilinear form on a finite-dimensional real vector space, let $x$ be a vector, and let $r\in\mathbb R$. The contravariant form $B^{-1}+r\,x\otimes x$ is nondegenerate if and only if $1+rB(x,x)\ne0$. When it is nondegenerate, its inverse $g$ satisfies
\begin{align}
g&=B-\frac{r}{1+rB(x,x)}B(x,\cdot)\otimes B(x,\cdot),\notag\\
g(x,\cdot)&=\frac{B(x,\cdot)}{1+rB(x,x)},\notag\\
g(x,x)&=\frac{B(x,x)}{1+rB(x,x)}.
\label{eq:rank-one-inverse}
\end{align}
The ordered signatures of $g$ and $B$ agree if and only if $1+rB(x,x)>0$. Moreover, any two nondegenerate symmetric forms satisfy
\begin{equation}
\operatorname{rank}(g^{-1}-B^{-1})
=\operatorname{rank}(g-B).
\label{eq:inverse-rank}
\end{equation}
\end{lemma}

\begin{proof}
Choose a basis and write the matrix of $B^{-1}+r\,x\otimes x$ as $(B^{ij}+rx^ix^j)$. In the multilinear column expansion of its determinant, a term using at least two columns from $(rx^ix^j)$ vanishes, since those columns are proportional to $x$. The remaining terms are
\begin{align*}
\det(B^{ij}+rx^ix^j)
&=\det(B^{ij})
+r\sum_{i=1}^n\sum_{j=1}^n
\operatorname{cof}_{ij}(B^{ab})x^ix^j\\
&=\det(B^{ij})
+r\det(B^{ij})\sum_{i=1}^n\sum_{j=1}^nB_{ji}x^ix^j\\
&=\det(B^{-1})\bigl(1+rB(x,x)\bigr).
\end{align*}
Here the cofactor identity follows from the inverse-matrix formula:
\[
\operatorname{cof}_{ij}(B^{ab})=\det(B^{ab})B_{ji}.
\]
Since $\det(B^{-1})\ne0$, this proves the nondegeneracy criterion.

Under this condition, multiply the proposed inverse by the original matrix. Every sum below runs over explicitly displayed indices:
\begin{align*}
&\sum_{k=1}^n(B^{ik}+rx^ix^k)
\left(B_{kj}-\frac{r}{1+rB(x,x)}
\sum_{a=1}^n\sum_{b=1}^n B_{ka}x^aB_{jb}x^b\right)\\
&=\sum_{k=1}^n B^{ik}B_{kj}
-\frac{r}{1+rB(x,x)}
\sum_{a=1}^n\sum_{b=1}^n
\left(\sum_{k=1}^n B^{ik}B_{ka}\right)x^aB_{jb}x^b\\
&\quad+rx^i\sum_{k=1}^n x^kB_{kj}
-\frac{r^2}{1+rB(x,x)}x^i
\sum_{b=1}^n\left(\sum_{k=1}^n\sum_{a=1}^n x^kB_{ka}x^a\right)B_{jb}x^b\\
&=\delta^i{}_j
-\frac{r}{1+rB(x,x)}x^i\sum_{b=1}^n B_{jb}x^b
+rx^i\sum_{b=1}^n B_{jb}x^b\\
&\quad-\frac{r^2B(x,x)}{1+rB(x,x)}x^i\sum_{b=1}^n B_{jb}x^b\\
&=\delta^i{}_j
+\frac{-r+r(1+rB(x,x))-r^2B(x,x)}{1+rB(x,x)}
x^i\sum_{b=1}^n B_{jb}x^b\\
&=\delta^i{}_j.
\end{align*}
This proves the inverse formula. Its contractions with $x$ give
\begin{align*}
g(x,\cdot)
&=B(x,\cdot)-\frac{rB(x,x)}{1+rB(x,x)}B(x,\cdot)\\
&=\frac{1+rB(x,x)-rB(x,x)}{1+rB(x,x)}B(x,\cdot)\\
&=\frac{B(x,\cdot)}{1+rB(x,x)},\\
g(x,x)
&=B(x,x)-\frac{rB(x,x)^2}{1+rB(x,x)}\\
&=\frac{B(x,x)(1+rB(x,x))-rB(x,x)^2}{1+rB(x,x)}\\
&=\frac{B(x,x)}{1+rB(x,x)}.
\end{align*}

If $1+rB(x,x)>0$, the path
$B^{-1}+tr\,x\otimes x$, $0\leq t\leq1$, remains nondegenerate, because
\[
1+trB(x,x)=(1-t)+t\bigl(1+rB(x,x)\bigr)>0.
\]
Indeed,
\[
\det(B^{-1}+tr\,x\otimes x)
=\det(B^{-1})\bigl(1+trB(x,x)\bigr)\ne0.
\]
If $\lambda_1(t),\ldots,\lambda_n(t)$ are the ordered eigenvalues of its real symmetric matrix, then
\[
\lambda_i(t)\ne0\quad(1\le i\le n),\qquad
\#\{i:\lambda_i(t)>0\}=\#\{i:\lambda_i(0)>0\},
\]
and the analogous equality holds for negative eigenvalues, by continuity. Thus the inertia is constant along the path. The eigenvalues of an invertible real symmetric matrix and its inverse are reciprocals, and
\[
\lambda>0\Longleftrightarrow\lambda^{-1}>0,\qquad
\lambda<0\Longleftrightarrow\lambda^{-1}<0.
\]
It follows that $g$ and $B$ have the same ordered signature.

If $1+rB(x,x)<0$, then $B(x,x)\ne0$ and
\[
V=\mathbb Rx\oplus\{u:B(x,u)=0\}
\]
is an orthogonal decomposition for both forms. To check the second assertion, for $u,v$ in the second summand we compute
\begin{align*}
g(u,v)
&=B(u,v)-\frac{r}{1+rB(x,x)}B(x,u)B(x,v)
=B(u,v),\\
g(x,u)
&=B(x,u)-\frac{rB(x,x)}{1+rB(x,x)}B(x,u)=0.
\end{align*}
Meanwhile,
\[
g(x,x)=\frac{B(x,x)}{1+rB(x,x)},\qquad
\operatorname{sgn}g(x,x)=-\operatorname{sgn}B(x,x).
\]
Thus, if the ordered signature of $B$ is $(p,q)$,
\begin{align*}
B(x,x)>0&\quad\Longrightarrow\quad\operatorname{sign}(g)=(p-1,q+1),\\
B(x,x)<0&\quad\Longrightarrow\quad\operatorname{sign}(g)=(p+1,q-1).
\end{align*}
Exactly one sign changes, so the ordered signatures differ.

Finally, as compositions of the associated linear maps,
\begin{align*}
-g^{-1}(g-B)B^{-1}
&=-g^{-1}gB^{-1}+g^{-1}BB^{-1}\\
&=-B^{-1}+g^{-1}=g^{-1}-B^{-1}.
\end{align*}
Multiplication by invertible maps preserves rank. Therefore
\begin{align*}
\operatorname{rank}(g^{-1}-B^{-1})
&=\operatorname{rank}\bigl(-g^{-1}(g-B)B^{-1}\bigr)\\
&=\operatorname{rank}(g-B),
\end{align*}
which proves \eqref{eq:inverse-rank}.
\end{proof}

\begin{theorem}[A metric correspondence on its exact domain]
\label{thm:metric-correspondence}
Fix smooth $X$, $c$, and $r$, and suppose that
\begin{equation}
\left(Xr+2cr+\frac2m\right)X\otimes X=0.
\label{eq:transport-tensor}
\end{equation}
Define
\begin{align*}
\mathcal C_c^r(X)
&=\{B\in\mathcal C_c(X):1+rB(X,X)\ne0\text{ everywhere}\},\\
\mathcal M_c^r(X;m)
&=\{g\in\mathcal M_c(X;m):1-rg(X,X)\ne0\text{ everywhere}\}.
\end{align*}
Then
\begin{equation}
T_r(B)^{-1}=B^{-1}+rX\otimes X
\label{eq:inverse-translation}
\end{equation}
defines a bijection $T_r:\mathcal C_c^r(X)\longrightarrow\mathcal M_c^r(X;m)$.
Its inverse is given by $B^{-1}=g^{-1}-rX\otimes X$. The covariant formula and the corresponding denominator identity are
\begin{align}
g&=B-\frac{r}{1+rB(X,X)}B(X,\cdot)\otimes B(X,\cdot),\notag\\
g(X,\cdot)&=\frac{B(X,\cdot)}{1+rB(X,X)},\notag\\
1-rg(X,X)&=\frac1{1+rB(X,X)}.
\label{eq:covariant-translation}
\end{align}
At each point, this correspondence preserves the ordered signature if and only if $1+rB(X,X)>0$.

Its restriction gives a bijection between all metrics in $\mathcal C_c(X)$ with $B(X,X)=0$ and all metrics in $\mathcal M_c(X;m)$ with $g(X,X)=0$. On these subsets, it preserves both the signature and the entire metric-dual one-form, and
\[
g=B-rB(X,\cdot)\otimes B(X,\cdot).
\]
\end{theorem}

\begin{proof}
Since $[X,X]=0$, the product rule gives
\begin{align*}
\mathcal L_X(X\otimes X)
&=[X,X]\otimes X+X\otimes[X,X]=0,\\
\mathcal L_X(rX\otimes X)
&=(Xr)X\otimes X+r[X,X]\otimes X+rX\otimes[X,X]\\
&=(Xr)X\otimes X.
\end{align*}
The cancellation can also be seen in coordinates:
\begin{align*}
\bigl(\mathcal L_X(rX\otimes X)\bigr)^{ij}
&=\sum_{k=1}^n X^k\partial_k(rX^iX^j)
-r\sum_{k=1}^n X^kX^j\partial_kX^i
-r\sum_{k=1}^n X^iX^k\partial_kX^j\\
&=\left(\sum_{k=1}^n X^k\partial_kr\right)X^iX^j
+r\sum_{k=1}^n X^k(\partial_kX^i)X^j
+r\sum_{k=1}^n X^kX^i\partial_kX^j\\
&\quad-r\sum_{k=1}^n X^kX^j\partial_kX^i
-r\sum_{k=1}^n X^iX^k\partial_kX^j\\
&=(Xr)X^iX^j.
\end{align*}
It follows from \eqref{eq:transport-tensor} that
\begin{align*}
\mathcal L_X(rX\otimes X)+2crX\otimes X
&=(Xr+2cr)X\otimes X\\
&=-\frac2mX\otimes X.
\end{align*}

For $B\in\mathcal C_c^r(X)$, define $g^{-1}$ by \eqref{eq:inverse-translation}. Its determinant is
\[
\det g^{-1}
=\det(B^{-1}+rX\otimes X)
=\det(B^{-1})\bigl(1+rB(X,X)\bigr)\ne0.
\]
Thus its inverse is a smooth metric. Substituting its definition into the inverse-metric equation gives
\begin{align*}
\mathcal L_Xg^{-1}+2cg^{-1}
&=\mathcal L_X(B^{-1}+rX\otimes X)
+2c(B^{-1}+rX\otimes X)\\
&=-2cB^{-1}+(Xr)X\otimes X
+2cB^{-1}+2crX\otimes X\\
&=(Xr+2cr)X\otimes X\\
&=-\frac2mX\otimes X.
\end{align*}
Lemma~\ref{lem:inverse-equation} now gives $g\in\mathcal M_c(X;m)$. Lemma~\ref{lem:rank-one-inversion} and contraction give
\begin{align*}
g&=B-\frac{r}{1+rB(X,X)}B(X,\cdot)\otimes B(X,\cdot),\\
g(X,\cdot)
&=B(X,\cdot)-\frac{rB(X,X)}{1+rB(X,X)}B(X,\cdot)\\
&=\frac{B(X,\cdot)}{1+rB(X,X)},\\
g(X,X)
&=B(X,X)-\frac{rB(X,X)^2}{1+rB(X,X)}\\
&=\frac{B(X,X)}{1+rB(X,X)},\\
1-rg(X,X)
&=1-\frac{rB(X,X)}{1+rB(X,X)}\\
&=\frac{1+rB(X,X)-rB(X,X)}{1+rB(X,X)}\\
&=\frac1{1+rB(X,X)}\ne0.
\end{align*}
This proves that $g$ belongs to the asserted target set.

For completeness, the same transformation can be checked directly in the covariant equation. The identities for an ordinary conformal metric give
\begin{align*}
\mathcal L_X\bigl(B(X,\cdot)\bigr)&=2cB(X,\cdot),\\
X\bigl(B(X,X)\bigr)&=2cB(X,X),\\
X\bigl(1+rB(X,X)\bigr)
&=(Xr)B(X,X)+rX\bigl(B(X,X)\bigr)\\
&=(Xr)B(X,X)+2crB(X,X).
\end{align*}
The quotient rule therefore gives
\begin{align*}
X\left(\frac{r}{1+rB(X,X)}\right)
&=\frac{(Xr)(1+rB(X,X))-rX(1+rB(X,X))}
{(1+rB(X,X))^2}\\
&=\frac{(Xr)(1+rB(X,X))
-r\bigl((Xr)B(X,X)+2crB(X,X)\bigr)}
{(1+rB(X,X))^2}\\
&=\frac{Xr+rB(X,X)Xr-rB(X,X)Xr-2cr^2B(X,X)}
{(1+rB(X,X))^2}\\
&=\frac{Xr-2cr^2B(X,X)}{(1+rB(X,X))^2}.
\end{align*}
Apply the product rule to the covariant formula:
\begin{align*}
\mathcal L_Xg
&=\mathcal L_XB
-X\left(\frac{r}{1+rB(X,X)}\right)B(X,\cdot)\otimes B(X,\cdot)\\
&\quad-\frac{r}{1+rB(X,X)}
\bigl(\mathcal L_X(B(X,\cdot))\otimes B(X,\cdot)
+B(X,\cdot)\otimes\mathcal L_X(B(X,\cdot))\bigr)\\
&=2cB-\frac{Xr-2cr^2B(X,X)}{(1+rB(X,X))^2}
B(X,\cdot)\otimes B(X,\cdot)\\
&\quad-\frac{4cr}{1+rB(X,X)}B(X,\cdot)\otimes B(X,\cdot)\\
&=2cB-\frac{Xr-2cr^2B(X,X)+4cr(1+rB(X,X))}
{(1+rB(X,X))^2}B(X,\cdot)\otimes B(X,\cdot)\\
&=2cB-\frac{Xr+4cr+2cr^2B(X,X)}{(1+rB(X,X))^2}
B(X,\cdot)\otimes B(X,\cdot).
\end{align*}
Substituting this expression and $g(X,\cdot)=B(X,\cdot)/(1+rB(X,X))$ into the modified conformal equation gives the residual
\begin{align*}
&\frac12\mathcal L_Xg-\frac1m g(X,\cdot)\otimes g(X,\cdot)-cg\\
&=cB-\frac{Xr+4cr+2cr^2B(X,X)}{2(1+rB(X,X))^2}
B(X,\cdot)\otimes B(X,\cdot)\\
&\quad-\frac{B(X,\cdot)\otimes B(X,\cdot)}{m(1+rB(X,X))^2}
-cB+\frac{cr}{1+rB(X,X)}B(X,\cdot)\otimes B(X,\cdot)\\
&=\frac{-Xr-4cr-2cr^2B(X,X)-2/m+2cr(1+rB(X,X))}
{2(1+rB(X,X))^2}B(X,\cdot)\otimes B(X,\cdot)\\
&=-\frac{Xr+2cr+2/m}{2(1+rB(X,X))^2}
B(X,\cdot)\otimes B(X,\cdot).
\end{align*}
This residual vanishes even at zeros of $X$. Indeed, its $(i,j)$ component is
\begin{align*}
&-\frac1{2(1+rB(X,X))^2}
\sum_{a=1}^n\sum_{b=1}^n B_{ia}B_{jb}
\left[\left(Xr+2cr+\frac2m\right)X^aX^b\right]\\
&\qquad=0
\end{align*}
by the tensor hypothesis \eqref{eq:transport-tensor}; no scalar cancellation at a zero of $X$ has been used.

Conversely, let $g\in\mathcal M_c^r(X;m)$ and set $B^{-1}=g^{-1}-rX\otimes X$. Its determinant is
\[
\det B^{-1}
=\det(g^{-1}-rX\otimes X)
=\det(g^{-1})\bigl(1-rg(X,X)\bigr)\ne0.
\]
It therefore defines a smooth metric $B$. Substitution gives
\begin{align*}
\mathcal L_XB^{-1}+2cB^{-1}
&=\mathcal L_X(g^{-1}-rX\otimes X)
+2c(g^{-1}-rX\otimes X)\\
&=\mathcal L_Xg^{-1}-(Xr)X\otimes X
+2cg^{-1}-2crX\otimes X\\
&=-\frac2mX\otimes X-(Xr)X\otimes X-2crX\otimes X\\
&=-\left(Xr+2cr+\frac2m\right)X\otimes X\\
&=0.
\end{align*}
Inverting this equation gives
\begin{align*}
(\mathcal L_XB)_{ij}
&=-\sum_{k=1}^n\sum_{\ell=1}^n
B_{ik}(\mathcal L_XB^{-1})^{k\ell}B_{\ell j}\\
&=-\sum_{k=1}^n\sum_{\ell=1}^n B_{ik}(-2cB^{k\ell})B_{\ell j}\\
&=2c\sum_{\ell=1}^n\delta_i{}^\ell B_{\ell j}=2cB_{ij}.
\end{align*}
Thus $B\in\mathcal C_c(X)$. Applying Lemma~\ref{lem:rank-one-inversion} with $g$ in place of $B$ and $-r$ in place of $r$ yields
\begin{align*}
B&=g+\frac{r}{1-rg(X,X)}g(X,\cdot)\otimes g(X,\cdot),\\
B(X,\cdot)
&=g(X,\cdot)+\frac{rg(X,X)}{1-rg(X,X)}g(X,\cdot)
=\frac{g(X,\cdot)}{1-rg(X,X)},\\
B(X,X)
&=g(X,X)+\frac{rg(X,X)^2}{1-rg(X,X)}
=\frac{g(X,X)}{1-rg(X,X)},\\
1+rB(X,X)
&=1+\frac{rg(X,X)}{1-rg(X,X)}\\
&=\frac{1-rg(X,X)+rg(X,X)}{1-rg(X,X)}
=\frac1{1-rg(X,X)}\ne0.
\end{align*}
Hence $B$ belongs to the stated source set.

The two inverse-metric translations cancel:
\[
(B^{-1}+rX\otimes X)-rX\otimes X=B^{-1},\qquad
(g^{-1}-rX\otimes X)+rX\otimes X=g^{-1}.
\]
The corresponding covariant compositions can be checked using the denominator identities:
\begin{align*}
&g+\frac{r}{1-rg(X,X)}g(X,\cdot)\otimes g(X,\cdot)\\
&=B-\frac{r}{1+rB(X,X)}B(X,\cdot)\otimes B(X,\cdot)\\
&\quad+\frac{r}{1/(1+rB(X,X))}
\frac{B(X,\cdot)\otimes B(X,\cdot)}{(1+rB(X,X))^2}\\
&=B-\frac{r}{1+rB(X,X)}B(X,\cdot)\otimes B(X,\cdot)
+\frac{r}{1+rB(X,X)}B(X,\cdot)\otimes B(X,\cdot)\\
&=B,
\end{align*}
and
\begin{align*}
&B-\frac{r}{1+rB(X,X)}B(X,\cdot)\otimes B(X,\cdot)\\
&=g+\frac{r}{1-rg(X,X)}g(X,\cdot)\otimes g(X,\cdot)\\
&\quad-\frac{r}{1/(1-rg(X,X))}
\frac{g(X,\cdot)\otimes g(X,\cdot)}{(1-rg(X,X))^2}\\
&=g+\frac{r}{1-rg(X,X)}g(X,\cdot)\otimes g(X,\cdot)
-\frac{r}{1-rg(X,X)}g(X,\cdot)\otimes g(X,\cdot)\\
&=g.
\end{align*}
Their signature criterion is exactly the one in Lemma~\ref{lem:rank-one-inversion}:
\[
\operatorname{sign}(g)=\operatorname{sign}(B)
\quad\Longleftrightarrow\quad 1+rB(X,X)>0.
\]

If $B(X,X)=0$, the forward formulas give
\begin{align*}
1+rB(X,X)&=1,\\
g&=B-rB(X,\cdot)\otimes B(X,\cdot),\\
g(X,\cdot)&=\frac{B(X,\cdot)}1=B(X,\cdot),\\
g(X,X)&=\frac01=0.
\end{align*}
Conversely, if $g(X,X)=0$, then
\begin{align*}
1-rg(X,X)&=1,\\
B&=g+rg(X,\cdot)\otimes g(X,\cdot),\\
B(X,\cdot)&=\frac{g(X,\cdot)}1=g(X,\cdot),\\
B(X,X)&=\frac01=0.
\end{align*}
Both denominators are positive, so neither null subset requires an additional nondegeneracy restriction, and both the signature and the dual one-form are preserved. The proof uses the tensor equation \eqref{eq:transport-tensor} throughout and remains valid at zeros of $X$.
\end{proof}

\begin{corollary}[A nonzero constant conformal factor]
\label{cor:constant-factor}
Suppose that $c\ne0$ is constant. There is a bijection
\begin{align*}
&\{B\in\mathcal C_c(X):B(X,X)\ne mc\text{ everywhere}\}\\
&\hspace{2em}\longleftrightarrow
\{g\in\mathcal M_c(X;m):g(X,X)\ne-mc\text{ everywhere}\}.
\end{align*}
The forward formulas are
\begin{align}
g^{-1}&=B^{-1}-\frac1{mc}X\otimes X,\notag\\
g&=B+\frac{B(X,\cdot)\otimes B(X,\cdot)}{mc-B(X,X)},
\label{eq:constant-factor-forward}
\end{align}
and the inverse formulas are
\begin{align}
B^{-1}&=g^{-1}+\frac1{mc}X\otimes X,\notag\\
B&=g-\frac{g(X,\cdot)\otimes g(X,\cdot)}{mc+g(X,X)}.
\label{eq:constant-factor-backward}
\end{align}
The signature is preserved precisely where $1-B(X,X)/(mc)>0$. On the null subsets, this is a global bijection that preserves the signature and the dual one-form, with
\begin{equation}
g=B+\frac1{mc}B(X,\cdot)\otimes B(X,\cdot).
\label{eq:constant-factor-null}
\end{equation}
Neither completeness of $X$ nor any additional invariance assumption is required.
\end{corollary}

\begin{proof}
Set $r=-1/(mc)$ in Theorem~\ref{thm:metric-correspondence}. Since $c$ is constant,
\begin{align*}
Xr+2cr+\frac2m
&=X\left(-\frac1{mc}\right)+2c\left(-\frac1{mc}\right)+\frac2m\\
&=0-\frac2m+\frac2m=0.
\end{align*}
The two nondegeneracy factors become
\begin{align*}
1+rB(X,X)
&=1-\frac{B(X,X)}{mc}=\frac{mc-B(X,X)}{mc},\\
1-rg(X,X)
&=1+\frac{g(X,X)}{mc}=\frac{mc+g(X,X)}{mc}.
\end{align*}
Since $mc\ne0$, their nonvanishing is equivalent to $B(X,X)\ne mc$ and $g(X,X)\ne-mc$, respectively.

Substituting $r=-1/(mc)$ into the forward formulas gives
\begin{align*}
g^{-1}
&=B^{-1}+\left(-\frac1{mc}\right)X\otimes X
=B^{-1}-\frac1{mc}X\otimes X,\\
g
&=B-\frac{-1/(mc)}{1-B(X,X)/(mc)}B(X,\cdot)\otimes B(X,\cdot)\\
&=B+\frac{1/(mc)}{(mc-B(X,X))/(mc)}B(X,\cdot)\otimes B(X,\cdot)\\
&=B+\frac{B(X,\cdot)\otimes B(X,\cdot)}{mc-B(X,X)}.
\end{align*}
The reverse formulas become
\begin{align*}
B^{-1}
&=g^{-1}-\left(-\frac1{mc}\right)X\otimes X
=g^{-1}+\frac1{mc}X\otimes X,\\
B
&=g+\frac{-1/(mc)}{1+g(X,X)/(mc)}g(X,\cdot)\otimes g(X,\cdot)\\
&=g-\frac{1/(mc)}{(mc+g(X,X))/(mc)}g(X,\cdot)\otimes g(X,\cdot)\\
&=g-\frac{g(X,\cdot)\otimes g(X,\cdot)}{mc+g(X,X)}.
\end{align*}
The signature condition is
\[
1+rB(X,X)>0
\quad\Longleftrightarrow\quad
1-\frac{B(X,X)}{mc}>0.
\]
For $B(X,X)=0$, the forward formula and its contraction reduce to
\begin{align*}
g&=B+\frac{B(X,\cdot)\otimes B(X,\cdot)}{mc-0}
=B+\frac1{mc}B(X,\cdot)\otimes B(X,\cdot),\\
g(X,\cdot)
&=B(X,\cdot)+\frac{B(X,X)}{mc}B(X,\cdot)
=B(X,\cdot),\\
g(X,X)&=B(X,X)+\frac{B(X,X)^2}{mc}=0.
\end{align*}
The inverse has the analogous nullity and dual-form identities, as already calculated in Theorem~\ref{thm:metric-correspondence}.

A nonzero constant conformal factor does not by itself imply $B(X,X)=0$. If $B(X,X)$ is also constant, then
\[
0=X\bigl(B(X,X)\bigr)=2cB(X,X),\qquad c\ne0
\quad\Longrightarrow\quad B(X,X)=0.
\]
The nullity assumption in \eqref{eq:constant-factor-null} must therefore be retained.

The exceptional value $g(X,X)=-mc$ makes the inverse construction genuinely degenerate, since
\begin{align*}
\det\left(g^{-1}+\frac1{mc}X\otimes X\right)
&=\det(g^{-1})\left(1+\frac{g(X,X)}{mc}\right)\\
&=\det(g^{-1})\left(1+\frac{-mc}{mc}\right)=0.
\end{align*}
There is also an explicit nonzero kernel covector. At such a point, $g(X,X)=-mc\ne0$ implies $X\ne0$, so $g(X,\cdot)\ne0$, and
\begin{align*}
\left(g^{-1}+\frac1{mc}X\otimes X\right)g(X,\cdot)
&=X+\frac{g(X,X)}{mc}X\\
&=X+\frac{-mc}{mc}X=0.
\end{align*}
\end{proof}

\begin{corollary}[Metric germs at a regular point]
\label{cor:regular-point-germs}
Let $p\in M$ satisfy $X_p\ne0$. Fix a germ of a smooth factor $c$ at $p$ and a local hypersurface $\Sigma$ through $p$ transverse to $X$. The initial-value problem
\begin{equation}
Xr+2cr=-\frac2m,\qquad r|_\Sigma=0
\label{eq:transport-initial-data}
\end{equation}
has a unique local smooth solution. For this fixed $r$, \eqref{eq:inverse-translation} gives a signature-preserving bijection between all ordinary conformal metric germs and all modified conformal metric germs with the same factor germ $c$. The corresponding representatives agree on $\Sigma$ as full covariant tensors, including their components in transverse directions.
\end{corollary}

\begin{proof}
Choose flowbox coordinates $(t,y)$ near $p$ such that $X=\partial_t$ and $\Sigma=\{t=0\}$. The initial-value problem becomes
\[
\partial_tr(t,y)+2c(t,y)r(t,y)=-\frac2m,
\qquad r(0,y)=0.
\]
The integrating factor satisfies
\[
\partial_t e^{2\int_0^t c(\sigma,y)\,d\sigma}
=2c(t,y)e^{2\int_0^t c(\sigma,y)\,d\sigma}.
\]
Multiplication of the equation by this factor gives
\begin{align*}
\partial_t\left(e^{2\int_0^t c(\sigma,y)\,d\sigma}r(t,y)\right)
&=e^{2\int_0^t c(\sigma,y)\,d\sigma}\partial_tr(t,y)\\
&\quad+2c(t,y)e^{2\int_0^t c(\sigma,y)\,d\sigma}r(t,y)\\
&=-\frac2m e^{2\int_0^t c(\sigma,y)\,d\sigma}.
\end{align*}
Integration from $0$ to $t$ yields
\begin{align*}
&e^{2\int_0^t c(\sigma,y)\,d\sigma}r(t,y)
-e^{2\int_0^0 c(\sigma,y)\,d\sigma}r(0,y)\\
&\qquad=-\frac2m\int_0^t e^{2\int_0^u c(\sigma,y)\,d\sigma}\,du.
\end{align*}
The initial term vanishes. Dividing by the positive integrating factor gives the only possible solution:
\[
r(t,y)=-\frac2m e^{-2\int_0^t c(\sigma,y)\,d\sigma}
\int_0^t e^{2\int_0^u c(\sigma,y)\,d\sigma}\,du.
\]
This expression is smooth. Its derivative is
\begin{align*}
\partial_tr(t,y)
&=-\frac2m\partial_t\left(
e^{-2\int_0^t c(\sigma,y)\,d\sigma}
\int_0^t e^{2\int_0^u c(\sigma,y)\,d\sigma}\,du\right)\\
&=\frac4m c(t,y)e^{-2\int_0^t c(\sigma,y)\,d\sigma}
\int_0^t e^{2\int_0^u c(\sigma,y)\,d\sigma}\,du\\
&\quad-\frac2m e^{-2\int_0^t c(\sigma,y)\,d\sigma}
e^{2\int_0^t c(\sigma,y)\,d\sigma}\\
&=\frac4m c(t,y)e^{-2\int_0^t c(\sigma,y)\,d\sigma}
\int_0^t e^{2\int_0^u c(\sigma,y)\,d\sigma}\,du-\frac2m.
\end{align*}
The other term in the equation is
\[
2c(t,y)r(t,y)
=-\frac4m c(t,y)e^{-2\int_0^t c(\sigma,y)\,d\sigma}
\int_0^t e^{2\int_0^u c(\sigma,y)\,d\sigma}\,du.
\]
Adding the two expressions gives
\begin{align*}
\partial_tr(t,y)+2c(t,y)r(t,y)
&=\frac4m c(t,y)e^{-2\int_0^t c(\sigma,y)\,d\sigma}
\int_0^t e^{2\int_0^u c(\sigma,y)\,d\sigma}\,du-\frac2m\\
&\quad-\frac4m c(t,y)e^{-2\int_0^t c(\sigma,y)\,d\sigma}
\int_0^t e^{2\int_0^u c(\sigma,y)\,d\sigma}\,du\\
&=-\frac2m.
\end{align*}
The initial value is
\[
r(0,y)=-\frac2m e^{-2\int_0^0 c(\sigma,y)\,d\sigma}
\int_0^0 e^{2\int_0^u c(\sigma,y)\,d\sigma}\,du=0.
\]
This proves existence. If $r_1$ and $r_2$ both solve the initial-value problem, their difference satisfies
\begin{align*}
\partial_t(r_1-r_2)+2c(r_1-r_2)&=0,\qquad (r_1-r_2)(0,y)=0,\\
\partial_t\left(e^{2\int_0^t c(\sigma,y)\,d\sigma}(r_1-r_2)\right)&=0,\\
e^{2\int_0^t c(\sigma,y)\,d\sigma}(r_1-r_2)&=0,\\
r_1-r_2&=0.
\end{align*}
Thus the local smooth solution is unique.

Substituting its explicit expression into the inverse-metric transformation gives
\begin{align*}
g^{-1}(t,y)
&=B^{-1}(t,y)\\
&\quad-\frac2m e^{-2\int_0^t c(\sigma,y)\,d\sigma}
\left(\int_0^t e^{2\int_0^u c(\sigma,y)\,d\sigma}\,du\right)
(X\otimes X)_{(t,y)}.
\end{align*}
The covariant expression is
\begin{align*}
g_{(t,y)}
&=B_{(t,y)}\\
&\quad+
\frac{\displaystyle\frac2m e^{-2\int_0^t c(\sigma,y)\,d\sigma}
\int_0^t e^{2\int_0^u c(\sigma,y)\,d\sigma}\,du}
{\displaystyle1-\frac2m e^{-2\int_0^t c(\sigma,y)\,d\sigma}
\left(\int_0^t e^{2\int_0^u c(\sigma,y)\,d\sigma}\,du\right)B_{(t,y)}(X,X)}\\
&\qquad\cdot B_{(t,y)}(X,\cdot)\otimes B_{(t,y)}(X,\cdot).
\end{align*}

For each ordinary conformal metric germ, the denominator at $p$ is
\[
1+r(p)B_p(X_p,X_p)
=1-\frac2m e^{-2\int_0^0 c(\sigma,y)\,d\sigma}
\left(\int_0^0 e^{2\int_0^u c(\sigma,y)\,d\sigma}\,du\right)B_p(X_p,X_p)
=1.
\]
It is positive after shrinking the representative's domain. In reverse,
\begin{align*}
B^{-1}(t,y)
&=g^{-1}(t,y)\\
&\quad+\frac2m e^{-2\int_0^t c(\sigma,y)\,d\sigma}
\left(\int_0^t e^{2\int_0^u c(\sigma,y)\,d\sigma}\,du\right)
(X\otimes X)_{(t,y)},\\
1-r(t,y)g_{(t,y)}(X,X)
&=1+\frac2m e^{-2\int_0^t c(\sigma,y)\,d\sigma}
\left(\int_0^t e^{2\int_0^u c(\sigma,y)\,d\sigma}\,du\right)
g_{(t,y)}(X,X).
\end{align*}
At $p$, this factor also equals $1-r(p)g_p(X_p,X_p)=1$, so it is positive on a sufficiently small neighborhood for each modified conformal metric germ. Theorem~\ref{thm:metric-correspondence} therefore applies in both directions and preserves the signature. The required neighborhood may depend on the chosen metric germ; no common neighborhood for all germs is asserted. Restriction to smaller neighborhoods commutes with both formulas, so the maps are well defined on germs.

Since $r|_\Sigma=0$, the covariant formula gives
\begin{align*}
g_{ij}|_\Sigma
&=B_{ij}|_\Sigma
-\left.\frac{r}{1+rB(X,X)}
\left(\sum_{a=1}^n B_{ia}X^a\right)
\left(\sum_{b=1}^n B_{jb}X^b\right)\right|_\Sigma\\
&=B_{ij}|_\Sigma.
\end{align*}
This holds for every pair of indices, including the transverse index. Thus the equality concerns sections of $S^2T^*M|_\Sigma$, rather than only their pullbacks to $T\Sigma$.
\end{proof}

\begin{remark}[Transport coefficients and zeros of the vector field]
\label{rem:transport-zeros}
On the regular set $\{X\ne0\}$, condition \eqref{eq:transport-tensor} is equivalent to the scalar equation $Xr+2cr=-2/m$. Indeed, at a regular point $p$ choose a covector $\lambda$ with $\lambda(X_p)\ne0$. Evaluation on $(\lambda,\lambda)$ gives
\begin{align*}
0
&=\left[\left(Xr+2cr+\frac2m\right)X\otimes X\right]_p(\lambda,\lambda)\\
&=\left(Xr(p)+2c(p)r(p)+\frac2m\right)\lambda(X_p)^2.
\end{align*}
Division by $\lambda(X_p)^2\ne0$ yields
\[
Xr(p)+2c(p)r(p)=-\frac2m.
\]
Conversely, multiplying this scalar equation by $X\otimes X$ gives \eqref{eq:transport-tensor}. At a zero of $X$, however, the condition reads
\[
\left(Xr(p)+2c(p)r(p)+\frac2m\right)0\otimes0=0.
\]
It is automatic and does not imply the scalar equation. Theorem~\ref{thm:metric-correspondence} assumes that a smooth $r$ satisfying its tensor condition has been supplied; it does not assert global existence of such a coefficient.

For example, suppose that $X_p=0$, $c(p)=0$, and regular points accumulate at $p$. If $p_j\to p$ is such a sequence and $r$ is a smooth solution of \eqref{eq:transport-tensor} near $p$, then
\[
dr_{p_j}(X_{p_j})+2c(p_j)r(p_j)+\frac2m=0.
\]
Passing to the limit gives
\begin{align*}
0
&=\lim_{j\to\infty}\left(dr_{p_j}(X_{p_j})+2c(p_j)r(p_j)+\frac2m\right)\\
&=dr_p(X_p)+2c(p)r(p)+\frac2m\\
&=dr_p(0)+2\cdot0\cdot r(p)+\frac2m
=\frac2m,
\end{align*}
a contradiction. This excludes the stated smooth inverse-metric translation near $p$; it does not exclude other metric constructions.

For a nonzero constant $c$, the coefficient $r=-1/(mc)$ is smooth everywhere. At a zero of $X$, the nondegeneracy factor in \eqref{eq:constant-factor-forward} is
\[
1-\frac{B_p(X_p,X_p)}{mc}
=1-\frac{B_p(0,0)}{mc}=1.
\]
By continuity, this factor remains positive on a sufficiently small neighborhood of that point, so the formula is smooth there and preserves the signature.

The covariant coefficient satisfies a Riccati equation. On the regular set, consider a smooth scalar $f$ and the formula
\[
g=B+fB(X,\cdot)\otimes B(X,\cdot),\qquad
1+fB(X,X)\ne0.
\]
Its determinant and contraction are
\begin{align*}
\det g&=\det B\bigl(1+fB(X,X)\bigr)\ne0,\\
g(X,\cdot)
&=B(X,\cdot)+fB(X,X)B(X,\cdot)\\
&=\bigl(1+fB(X,X)\bigr)B(X,\cdot).
\end{align*}
Using the ordinary conformal identities, its Lie derivative is
\begin{align*}
\mathcal L_Xg
&=\mathcal L_XB+(Xf)B(X,\cdot)\otimes B(X,\cdot)\\
&\quad+f\bigl(\mathcal L_X(B(X,\cdot))\otimes B(X,\cdot)
+B(X,\cdot)\otimes\mathcal L_X(B(X,\cdot))\bigr)\\
&=2cB+(Xf)B(X,\cdot)\otimes B(X,\cdot)
+4cfB(X,\cdot)\otimes B(X,\cdot)\\
&=2cB+(Xf+4cf)B(X,\cdot)\otimes B(X,\cdot).
\end{align*}
The right-hand side of the modified conformal equation, multiplied by $2$, is
\begin{align*}
2cg+\frac2m g(X,\cdot)\otimes g(X,\cdot)
&=2cB+2cfB(X,\cdot)\otimes B(X,\cdot)\\
&\quad+\frac2m\bigl(1+fB(X,X)\bigr)^2
B(X,\cdot)\otimes B(X,\cdot).
\end{align*}
Their difference is therefore
\begin{align*}
&\mathcal L_Xg-2cg-\frac2m g(X,\cdot)\otimes g(X,\cdot)\\
&\qquad=\left[Xf+2cf-\frac2m\bigl(1+fB(X,X)\bigr)^2\right]
B(X,\cdot)\otimes B(X,\cdot).
\end{align*}
At a regular point choose $u$ with $B(X,u)\ne0$. Evaluation on $(u,u)$ and division by $B(X,u)^2$ show that this tensor vanishes if and only if
\[
Xf+2cf=\frac2m\bigl(1+fB(X,X)\bigr)^2.
\]
For the particular coefficient $f=-r/(1+rB(X,X))$, the identity is
\begin{align}
&X\left(-\frac{r}{1+rB(X,X)}\right)
+2c\left(-\frac{r}{1+rB(X,X)}\right)\notag\\
&\qquad=\frac2m\left(1-\frac{rB(X,X)}{1+rB(X,X)}\right)^2,
\qquad 1-\frac{rB(X,X)}{1+rB(X,X)}\ne0.
\label{eq:covariant-riccati}
\end{align}
The substitution can be expanded explicitly. First,
\begin{align*}
X\bigl(1+rB(X,X)\bigr)
&=(Xr)B(X,X)+rX\bigl(B(X,X)\bigr)\\
&=(Xr)B(X,X)+2crB(X,X).
\end{align*}
Thus the quotient rule gives
\begin{align*}
X\left(-\frac{r}{1+rB(X,X)}\right)
&=-\frac{Xr}{1+rB(X,X)}
+\frac{rX(1+rB(X,X))}{(1+rB(X,X))^2}\\
&=-\frac{Xr}{1+rB(X,X)}
+\frac{r\bigl((Xr)B(X,X)+2crB(X,X)\bigr)}{(1+rB(X,X))^2}\\
&=\frac{-(Xr)(1+rB(X,X))+r(Xr)B(X,X)+2cr^2B(X,X)}
{(1+rB(X,X))^2}\\
&=\frac{-Xr+2cr^2B(X,X)}{(1+rB(X,X))^2}.
\end{align*}
Adding the linear term gives
\begin{align*}
X\left(-\frac{r}{1+rB(X,X)}\right)
-\frac{2cr}{1+rB(X,X)}
&=\frac{-Xr+2cr^2B(X,X)}{(1+rB(X,X))^2}
-\frac{2cr}{1+rB(X,X)}\\
&=\frac{-Xr+2cr^2B(X,X)-2cr(1+rB(X,X))}{(1+rB(X,X))^2}\\
&=\frac{-Xr-2cr}{(1+rB(X,X))^2}.
\end{align*}
On the right-hand side,
\begin{align*}
1-\frac{rB(X,X)}{1+rB(X,X)}
&=\frac{1+rB(X,X)-rB(X,X)}{1+rB(X,X)}\\
&=\frac1{1+rB(X,X)},\\
\frac2m\left(1-\frac{rB(X,X)}{1+rB(X,X)}\right)^2
&=\frac{2/m}{(1+rB(X,X))^2}.
\end{align*}
Consequently, \eqref{eq:covariant-riccati} is equivalent to
\begin{align*}
\frac{-Xr-2cr}{(1+rB(X,X))^2}
&=\frac{2/m}{(1+rB(X,X))^2},\\
-Xr-2cr&=\frac2m,\\
Xr+2cr&=-\frac2m.
\end{align*}
The denominator is nonzero throughout this calculation. Thus the fractional change of coefficient converts the covariant Riccati equation into a linear transport equation.

Finally, if $c$ is smooth and nowhere zero but is not assumed constant, differentiation gives
\begin{align*}
X\left(-\frac1{mc}\right)
&=-\frac1mX(c^{-1})
=-\frac1m\left(-\frac{Xc}{c^2}\right)
=\frac{Xc}{mc^2},\\
X\left(-\frac1{mc}\right)
+2c\left(-\frac1{mc}\right)+\frac2m
&=\frac{Xc}{mc^2}-\frac2m+\frac2m
=\frac{Xc}{mc^2}.
\end{align*}
Hence $r=-1/(mc)$ satisfies the scalar transport equation precisely when $Xc=0$. In the tensor condition \eqref{eq:transport-tensor}, the residual is
\[
\left[X\left(-\frac1{mc}\right)
+2c\left(-\frac1{mc}\right)+\frac2m\right]X\otimes X
=\frac{Xc}{mc^2}X\otimes X.
\]
The derivative of $c$ along $X$ cannot in general be omitted.
\end{remark}

\section{Rank-one deformations of Killing metrics}
\label{sec:rank-one}

Throughout this section, \(X\) is nowhere zero and \(B\) is a metric
such that \(\mathcal L_XB=0\). We impose a rank bound on \(g-B\)
without prescribing its direction.

\begin{theorem}[Classification under a rank-one bound]
\label{thm:rank-one-classification}
Suppose that \(n\geq3\). A metric \(g\) and a smooth function \(c\)
satisfy
\[
 \frac12\mathcal L_Xg-\frac1m
 g(X,\cdot)\otimes g(X,\cdot)=cg,
 \qquad \operatorname{rank}(g-B)\leq1
\]
everywhere if and only if \(c=0\) and there is a smooth function \(r\)
such that
\begin{equation}
 Xr=-\frac2m,\qquad 1+rB(X,X)\ne0,\qquad
 g^{-1}=B^{-1}+rX\otimes X.
 \label{eq:rank-one-classification}
\end{equation}
The function \(r\) is unique, and
\begin{equation}
 g=B-\frac{r}{1+rB(X,X)}
 B(X,\cdot)\otimes B(X,\cdot).
 \label{eq:rank-one-covariant}
\end{equation}
If \(c=0\) is assumed from the outset, the same classification holds
in dimension \(2\).
\end{theorem}

\begin{proof}
We first show that the rank bound forces \(c=0\) when \(n\geq3\).
Set \(K=g-B\). The set on which \(K\) has rank one is open.
It is also dense. Otherwise, \(g=B\) on a nonempty open set, where
the modified conformal equation would give
\begin{align*}
0&=\frac12\mathcal L_XB-\frac1m
B(X,\cdot)\otimes B(X,\cdot)-cB\\
&=-\frac1m B(X,\cdot)\otimes B(X,\cdot)-cB.
\end{align*}
Thus \(cB=-B(X,\cdot)\otimes B(X,\cdot)/m\).
The right-hand side has rank one, since \(X\ne0\) and \(B\) is
nondegenerate. The left-hand side has rank zero or \(n\), which
is impossible.

Near any rank-one point, choose a smooth vector field \(W\) such
that \(K(W,W)\ne0\), and write
\[
 \alpha=K(W,\cdot),\qquad \sigma=K(W,W).
\]
Such a \(W\) exists by polarization: if \(K(W,W)=0\) for all
\(W\), then
\[
2K(U,V)=K(U+V,U+V)-K(U,U)-K(V,V)=0,
\]
contrary to \(K\ne0\).
The vanishing of the two-by-two minors gives
\begin{align*}
0&=\det\begin{pmatrix}
 K(U,V)&K(U,W)\\ K(W,V)&K(W,W)
 \end{pmatrix}\\
 &=\sigma K(U,V)-\alpha(U)\alpha(V).
\end{align*}
Therefore
\[
 K(U,V)=\frac{\alpha(U)\alpha(V)}{\sigma},
 \qquad K=\frac1{\sigma}\alpha\otimes\alpha.
\]
Let \(H=\ker\alpha\), and take smooth local sections \(U,V\) of
\(H\). The full product rule gives
\begin{align*}
(\mathcal L_XK)(U,V)
={}&X(\sigma^{-1})\alpha(U)\alpha(V)\\
 &+\sigma^{-1}(\mathcal L_X\alpha)(U)\alpha(V)
 +\sigma^{-1}\alpha(U)(\mathcal L_X\alpha)(V)\\
={}&X(\sigma^{-1})\,0\cdot0
 +\sigma^{-1}(\mathcal L_X\alpha)(U)\,0
 +\sigma^{-1}0\,(\mathcal L_X\alpha)(V)=0.
\end{align*}
On the other hand,
\begin{align*}
\frac12(\mathcal L_XK)(U,V)
 &=cg(U,V)+\frac1m g(X,U)g(X,V)\\
 &=c\left(B(U,V)+\frac{\alpha(U)\alpha(V)}{\sigma}\right)
   +\frac1m g(X,U)g(X,V)\\
 &=cB(U,V)+\frac1m g(X,U)g(X,V).
\end{align*}
Combining these calculations gives
\begin{equation}
 cB(U,V)=-\frac1m g(X,U)g(X,V),\qquad U,V\in H.
 \label{eq:rank-one-hyperplane}
\end{equation}
Here \(H\) has dimension \(n-1\), and
\[
 H^{\perp_B}=\mathbb R B^\sharp\alpha.
\]
Indeed,
\[
B(B^\sharp\alpha,U)=\alpha(U)=0\quad(U\in H),
\qquad \dim H^{\perp_B}=n-\dim H=1.
\]
Consequently, the radical \(H\cap H^{\perp_B}\) has dimension at
most one, and
\[
\operatorname{rank}(B|_H)
=\dim H-\dim(H\cap H^{\perp_B})\geq(n-1)-1=n-2.
\]
For \(n\geq4\), equation~\eqref{eq:rank-one-hyperplane} rules out
\(c\ne0\), because its left-hand side would have rank at least
two and its right-hand side has rank at most one.

Suppose that \(n=3\), and assume that \(c\ne0\) at a rank-one
point. We may work on a neighborhood where both conditions
persist. Equation~\eqref{eq:rank-one-hyperplane} forces
\(\operatorname{rank}(B|_H)=1\), because
\[
1=n-2\leq\operatorname{rank}(B|_H)\leq1.
\]
The radical therefore has dimension \(2-1=1\), and coincides
with \(H^{\perp_B}\). Thus, with \(v=B^\sharp\alpha\),
\[
 v\in H,\qquad \alpha(v)=B^{-1}(\alpha,\alpha)=0.
\]
Evaluating~\eqref{eq:rank-one-hyperplane} on \((v,v)\) gives
\[
0=cB(v,v)=c\alpha(v)=-\frac1m g(X,v)^2,
\qquad g(X,v)=0.
\]
On the other hand,
\begin{align*}
g(X,v)&=B(X,v)+\frac{\alpha(X)\alpha(v)}{\sigma}\\
&=\alpha(X)+\frac{\alpha(X)\cdot0}{\sigma}
=\alpha(X).
\end{align*}
It follows that \(\alpha(X)=0\), and hence
\[
g(X,U)=B(X,U)+\frac{\alpha(X)\alpha(U)}{\sigma}
=B(X,U)+\frac{0\cdot\alpha(U)}{\sigma}=B(X,U)
\]
for every \(U\). Thus
\begin{equation}
 g(X,\cdot)=B(X,\cdot).
 \label{eq:three-dimensional-dual}
\end{equation}
To check the metric densities explicitly, write
\(g_{ij}=B_{ij}+\alpha_i\alpha_j/\sigma\).
Multilinearity of the determinant gives
\begin{align*}
\det(g_{ij})
 &=\det(B_{ij})+\frac1{\sigma}
 \sum_{i,j=1}^3\operatorname{cof}_{ij}(B)\alpha_i\alpha_j\\
 &=\det(B_{ij})\left(
 1+\frac1{\sigma}\sum_{i,j=1}^3B^{ji}\alpha_i\alpha_j\right)\\
 &=\det(B_{ij})\left(1+
 \frac{B^{-1}(\alpha,\alpha)}{\sigma}\right)
 =\det(B_{ij}).
\end{align*}
Terms involving two or more changed columns vanish because
those columns are proportional to the column \((\alpha_i)\).
Here \(\operatorname{cof}_{ij}(B)=\det(B)B^{ji}\).
Therefore
\[
\sqrt{|\det(g_{ij})|}=\sqrt{|\det(B_{ij})|}.
\]
Using the coordinate formula for divergence now gives
\begin{align*}
\operatorname{div}_gX
 &=\frac1{\sqrt{|\det g|}}\sum_{i=1}^3
 \partial_i\left(\sqrt{|\det g|}X^i\right)\\
 &=\frac1{\sqrt{|\det B|}}\sum_{i=1}^3
 \partial_i\left(\sqrt{|\det B|}X^i\right)
 =\operatorname{div}_BX\\
 &=\frac12\sum_{i,j=1}^3B^{ij}(\mathcal L_XB)_{ij}=0.
\end{align*}
The trace of the modified conformal equation is
\begin{align*}
0=\operatorname{div}_gX
 &=\frac12\sum_{i,j=1}^3g^{ij}(\mathcal L_Xg)_{ij}\\
 &=c\sum_{i,j=1}^3g^{ij}g_{ij}
 +\frac1m\sum_{i,j,a,b=1}^3g^{ij}g_{ia}g_{jb}X^aX^b\\
 &=3c+\frac1m\sum_{j,a,b=1}^3\delta_a^j g_{jb}X^aX^b\\
 &=3c+\frac1m\sum_{a,b=1}^3g_{ab}X^aX^b
 =3c+\frac{g(X,X)}m.
\end{align*}
Meanwhile, \eqref{eq:three-dimensional-dual} and
\(\mathcal L_XB=0\) imply
\begin{align*}
\mathcal L_X(g(X,\cdot))
 &=(\mathcal L_Xg)(X,\cdot)+g([X,X],\cdot)\\
 &=\left(2cg+\frac2m g(X,\cdot)\otimes g(X,\cdot)\right)
 (X,\cdot)\\
 &=2cg(X,\cdot)+\frac2m g(X,X)g(X,\cdot)\\
 &=2\left(c+\frac{g(X,X)}m\right)g(X,\cdot),
\end{align*}
whereas
\[
\mathcal L_X(g(X,\cdot))
=\mathcal L_X(B(X,\cdot))
=(\mathcal L_XB)(X,\cdot)+B([X,X],\cdot)=0.
\]
The one-form \(g(X,\cdot)\) is nonzero. Hence
\(c+g(X,X)/m=0\), and subtracting the two scalar identities gives
\[
0=\left(3c+\frac{g(X,X)}m\right)
 -\left(c+\frac{g(X,X)}m\right)=2c,
\]
a contradiction. Therefore \(c=0\) on the rank-one set
in every dimension \(n\geq3\). Density and continuity give
\(c=0\) on \(M\).

We now assume \(c=0\); the remainder of the argument also applies
when \(n=2\). By Lemma~\ref{lem:inverse-equation}, the tensor
\[
 T=g^{-1}-B^{-1}
\]
satisfies
\begin{equation}
 \mathcal L_XT=-\frac2mX\otimes X,\qquad
 \operatorname{rank}T\leq1.
 \label{eq:rank-one-inverse-evolution}
\end{equation}
The rank bound follows from
\begin{align*}
-g^{-1}(g-B)B^{-1}
 &=-g^{-1}gB^{-1}+g^{-1}BB^{-1}\\
 &=-B^{-1}+g^{-1}=g^{-1}-B^{-1}.
\end{align*}
On a small neighborhood where \(T\) has rank one, it has the form
\(T=fV\otimes V\), where \(V\) is a nonzero smooth vector field
and \(f\) is smooth and nowhere zero. More explicitly, choose
a smooth one-form \(\eta\) with \(T(\eta,\eta)\ne0\), and set
\[
V=T^\sharp\eta,\qquad f=\frac1{T(\eta,\eta)}.
\]
Here \(T^\sharp\eta\) means contraction of the displayed tensor:
\[
V=\sum_{i,j=1}^n(g^{ij}-B^{ij})\eta_j\partial_{x^i}.
\]
The same minor identity gives, for covectors \(\lambda,\mu\),
\begin{align*}
0&=\det\begin{pmatrix}
 T(\lambda,\mu)&T(\lambda,\eta)\\
 T(\eta,\mu)&T(\eta,\eta)
 \end{pmatrix},\\
T(\lambda,\mu)
 &=\frac{T(\lambda,\eta)T(\eta,\mu)}{T(\eta,\eta)}
 =f\,\lambda(V)\mu(V).
\end{align*}
If a covector \(\lambda\) annihilates \(V\), the full product rule
gives
\begin{align*}
(\mathcal L_XT)(\lambda,\lambda)
={}&(Xf)\lambda(V)^2\\
 &+f\,\lambda([X,V])\lambda(V)
 +f\,\lambda(V)\lambda([X,V])\\
={}&(Xf)\,0^2+f\,\lambda([X,V])\,0
 +f\,0\,\lambda([X,V])=0.
\end{align*}
Equation~\eqref{eq:rank-one-inverse-evolution} then gives
\[
0=-\frac2m(X\otimes X)(\lambda,\lambda)
 =-\frac2m\lambda(X)^2,\qquad \lambda(X)=0.
\]
Thus every covector that annihilates \(V\)
also annihilates \(X\), so \(X\in\mathbb RV\).
It follows that \(T\) is a scalar multiple of \(X\otimes X\)
on this neighborhood.

To obtain a global coefficient that is smooth even where the rank
drops, choose an auxiliary Riemannian metric \(h\), and define
\[
 \theta=\frac{h(X,\cdot)}{h(X,X)},\qquad
 r=T(\theta,\theta).
\]
The denominator is positive, and
\[
\theta(X)=\frac{h(X,X)}{h(X,X)}=1,\qquad
r=\frac{g^{-1}(h(X,\cdot),h(X,\cdot))
 -B^{-1}(h(X,\cdot),h(X,\cdot))}{h(X,X)^2}.
\]
In particular, \(r\) is globally smooth.
On a rank-one neighborhood, write
\[
V=aX,\qquad a=\frac{h(X,V)}{h(X,X)}.
\]
Then
\begin{align*}
T&=f(aX)\otimes(aX)=fa^2X\otimes X,\\
r=T(\theta,\theta)
 &=fa^2\theta(X)\theta(X)=fa^2.
\end{align*}
Where \(T=0\), the same definition gives \(r=0\).
Hence \(T=rX\otimes X\) everywhere.
If also \(T=\widetilde rX\otimes X\), then
\[
0=((r-\widetilde r)X\otimes X)(\theta,\theta)
 =(r-\widetilde r)\theta(X)^2=r-\widetilde r,
\]
so \(r\) is unique. Substituting this expression into
\eqref{eq:rank-one-inverse-evolution}, and using \([X,X]=0\), gives
\begin{align*}
-\frac2mX\otimes X
 &=\mathcal L_X(rX\otimes X)\\
 &=(Xr)X\otimes X+r[X,X]\otimes X+rX\otimes[X,X]\\
 &=(Xr)X\otimes X.
\end{align*}
Contraction with \((\theta,\theta)\) yields
\[
-\frac2m=-\frac2m\theta(X)^2=(Xr)\theta(X)^2=Xr.
\]
Nondegeneracy and the
rank-one inversion formula from Section~2 give
\(1+rB(X,X)\ne0\) and~\eqref{eq:rank-one-covariant}.

Conversely, suppose that \(r\) satisfies
\eqref{eq:rank-one-classification}. The stated nondegeneracy
condition defines a metric \(g\), and
\[
 \mathcal L_Xg^{-1}
 =\mathcal L_XB^{-1}+(Xr)X\otimes X
 =-\frac2mX\otimes X.
\]
Lemma~\ref{lem:inverse-equation} gives the modified conformal
equation with factor zero.
For a covariant verification, the Killing condition gives
\begin{align*}
X(B(X,X))&=(\mathcal L_XB)(X,X)+2B([X,X],X)=0,\\
\mathcal L_X(B(X,\cdot))
 &=(\mathcal L_XB)(X,\cdot)+B([X,X],\cdot)=0.
\end{align*}
The quotient rule gives
\begin{align*}
X\left(\frac{r}{1+rB(X,X)}\right)
 &=\frac{(Xr)(1+rB(X,X))
 -r\bigl((Xr)B(X,X)+rX(B(X,X))\bigr)}
 {(1+rB(X,X))^2}\\
 &=\frac{(-2/m)(1+rB(X,X))-r(-2/m)B(X,X)}
 {(1+rB(X,X))^2}\\
 &=-\frac2{m(1+rB(X,X))^2}.
\end{align*}
Differentiating~\eqref{eq:rank-one-covariant} therefore gives
\begin{align*}
\mathcal L_Xg
={}&\mathcal L_XB-
X\left(\frac{r}{1+rB(X,X)}\right)
B(X,\cdot)\otimes B(X,\cdot)\\
 &-\frac{r}{1+rB(X,X)}
 \bigl(\mathcal L_X(B(X,\cdot))\otimes B(X,\cdot)
 +B(X,\cdot)\otimes\mathcal L_X(B(X,\cdot))\bigr)\\
={}&\frac2{m(1+rB(X,X))^2}
B(X,\cdot)\otimes B(X,\cdot).
\end{align*}
The contraction is
\begin{align*}
g(X,\cdot)
 &=B(X,\cdot)-\frac{rB(X,X)}{1+rB(X,X)}B(X,\cdot)\\
 &=\frac{B(X,\cdot)}{1+rB(X,X)}.
\end{align*}
Consequently,
\begin{align*}
\frac12\mathcal L_Xg-\frac1m g(X,\cdot)\otimes g(X,\cdot)
 &=
\frac{B(X,\cdot)\otimes B(X,\cdot)}
 {m(1+rB(X,X))^2}\\
 &\quad-\frac1m
 \frac{B(X,\cdot)}{1+rB(X,X)}
 \otimes\frac{B(X,\cdot)}{1+rB(X,X)}=0.
\end{align*}
Finally,
\[
g-B=-\frac{r}{1+rB(X,X)}
B(X,\cdot)\otimes B(X,\cdot)
\]
has rank at most one. This completes the converse.
\end{proof}

\begin{theorem}[Global criterion for a complete flow]
\label{thm:complete-rank-one}
Assume in addition that \(X\) is complete. If \(n\geq3\), there
exist a global metric \(g\) and a smooth function \(c\) satisfying
\[
 \frac12\mathcal L_Xg-\frac1m
 g(X,\cdot)\otimes g(X,\cdot)=cg,\qquad
 \operatorname{rank}(g-B)\leq1
\]
if and only if
\begin{equation}
 B(X,X)=0,\qquad
 \text{there is }\tau\in C^\infty(M)\text{ with }X\tau=1.
 \label{eq:complete-rank-one-criterion}
\end{equation}
In dimension \(2\), the same assertion holds for solutions whose
factor is required to be zero.
In these respective ranges, every solution has \(c=0\), and all
the metrics are parametrized bijectively by the functions \(\tau\)
with \(X\tau=1\), through
\begin{equation}
 g=B+\frac{2\tau}{m}B(X,\cdot)\otimes B(X,\cdot).
 \label{eq:complete-rank-one-metric}
\end{equation}
Each such \(g\) has the same signature as \(B\), and
\(g(X,\cdot)=B(X,\cdot)\).

The existence of a smooth function \(\tau\) with \(X\tau=1\)
is equivalent to the existence of a smooth manifold \(\Sigma\)
and a diffeomorphism
\[
 F:\mathbb R\times\Sigma\longrightarrow M,\qquad
 F_*(\partial_t)=X.
\]
\end{theorem}

\begin{proof}
Let \(\Phi_t\) be the flow of \(X\), defined for all
\(t\in\mathbb R\). By
Theorem~\ref{thm:rank-one-classification},
\begin{align*}
\frac d{dt}r(\Phi_t(p))&=(Xr)(\Phi_t(p))=-\frac2m,\\
r(\Phi_t(p))-r(p)&=\int_0^t-\frac2m\,ds=-\frac{2t}m,\\
r(\Phi_t(p))&=r(p)-\frac{2t}{m}.
\end{align*}
Since \(X\) is Killing for \(B\),
\begin{align*}
\frac d{dt}(B(X,X)(\Phi_t(p)))
 &=X(B(X,X))(\Phi_t(p))\\
 &=\bigl((\mathcal L_XB)(X,X)+2B([X,X],X)\bigr)(\Phi_t(p))
 =0,
\end{align*}
and therefore
\[
B(X,X)(\Phi_t(p))=B(X,X)(p).
\]
Thus
\[
 1+r(\Phi_t(p))B(X,X)(\Phi_t(p))
 =1+\left(r(p)-\frac{2t}{m}\right)B(X,X)(p).
\]
If \(B(X,X)(p)\ne0\), this expression vanishes at
\[
 t=\frac m2\left(r(p)+\frac1{B(X,X)(p)}\right),
\]
because substitution at this time gives
\begin{align*}
&1+\left[r(p)-\frac2m\frac m2
 \left(r(p)+\frac1{B(X,X)(p)}\right)\right]B(X,X)(p)\\
&\qquad=1+r(p)B(X,X)(p)-r(p)B(X,X)(p)-1=0.
\end{align*}
This contradicts nondegeneracy. Therefore \(B(X,X)=0\) everywhere.
Setting \(\tau=-mr/2\) gives
\[
X\tau=-\frac m2Xr=-\frac m2\left(-\frac2m\right)=1.
\]
Substitution into~\eqref{eq:rank-one-covariant} gives
\begin{align*}
g&=B-\frac{-2\tau/m}{1+(-2\tau/m)\cdot0}
B(X,\cdot)\otimes B(X,\cdot)\\
&=B+\frac{2\tau}mB(X,\cdot)\otimes B(X,\cdot),
\end{align*}
which is~\eqref{eq:complete-rank-one-metric}.

Conversely, if \eqref{eq:complete-rank-one-criterion} holds, then
\(r=-2\tau/m\) satisfies
\[
Xr=-\frac2mX\tau=-\frac2m,\qquad
1+rB(X,X)=1-\frac{2\tau}m\cdot0=1.
\]
Theorem~\ref{thm:rank-one-classification} gives the required metric.
The rank-one signature criterion shows that its signature agrees
with that of \(B\). The contraction and derivative can also be
checked explicitly:
\begin{align*}
g(X,\cdot)
 &=B(X,\cdot)+\frac{2\tau}mB(X,X)B(X,\cdot)
 =B(X,\cdot),\\
\mathcal L_Xg
 &=\mathcal L_XB+\frac2m(X\tau)B(X,\cdot)\otimes B(X,\cdot)\\
 &\quad+\frac{2\tau}m\bigl(
 \mathcal L_X(B(X,\cdot))\otimes B(X,\cdot)
 +B(X,\cdot)\otimes\mathcal L_X(B(X,\cdot))\bigr)\\
 &=0+\frac2m\cdot1\,B(X,\cdot)\otimes B(X,\cdot)
 +\frac{2\tau}m(0\otimes B(X,\cdot)+B(X,\cdot)\otimes0)\\
 &=\frac2m g(X,\cdot)\otimes g(X,\cdot).
\end{align*}
If two functions \(\tau_1,\tau_2\) yield the same metric, then
\[
\frac{2(\tau_1-\tau_2)}m B(X,\cdot)\otimes B(X,\cdot)=0.
\]
Choose at the point in question a vector \(W\) with \(B(X,W)\ne0\).
Evaluation gives
\[
\frac{2(\tau_1-\tau_2)}m B(X,W)^2=0,\qquad \tau_1=\tau_2.
\]
This proves the asserted parametrization.

It remains to prove the product characterization. Integration of
\(X\tau=1\) gives
\begin{align*}
\frac d{dt}\tau(\Phi_t(p))&=(X\tau)(\Phi_t(p))=1,\\
\tau(\Phi_t(p))-\tau(p)&=\int_0^t1\,ds=t,\\
\tau(\Phi_t(p))&=\tau(p)+t.
\end{align*}
The level set \(\Sigma=\tau^{-1}(0)\) is an embedded hypersurface,
because \(d\tau(X)=1\). Moreover,
\[
\tau(\Phi_{-\tau(p)}(p))=\tau(p)-\tau(p)=0.
\]
An orbit reaches \(\Sigma\) at a time \(t\) only if
\[
0=\tau(\Phi_t(p))=\tau(p)+t,\qquad t=-\tau(p),
\]
so each orbit meets \(\Sigma\) exactly once.
The map
\[
 F(t,q)=\Phi_t(q)
\]
has the smooth inverse
\[
 F^{-1}(p)=\bigl(\tau(p),\Phi_{-\tau(p)}(p)\bigr).
\]
Indeed, substitution in one direction gives
\begin{align*}
F\bigl(\tau(p),\Phi_{-\tau(p)}(p)\bigr)
 &=\Phi_{\tau(p)}(\Phi_{-\tau(p)}(p))\\
 &=\Phi_{\tau(p)-\tau(p)}(p)=p.
\end{align*}
In the other direction, \(q\in\Sigma\) implies
\begin{align*}
\tau(F(t,q))&=\tau(\Phi_t(q))=\tau(q)+t=t,\\
\Phi_{-\tau(F(t,q))}(F(t,q))
 &=\Phi_{-t}(\Phi_t(q))=\Phi_{-t+t}(q)=q.
\end{align*}
Both maps are smooth. Their action on the vector field is
\[
F_*(\partial_t)_{(t,q)}
=\frac d{ds}\bigg|_{s=0}F(t+s,q)
=\frac d{ds}\bigg|_{s=0}\Phi_{t+s}(q)
=X_{\Phi_t(q)}.
\]
Conversely, a diffeomorphism with this property gives
\(\tau=\operatorname{pr}_{\mathbb R}\circ F^{-1}\), for which
\begin{align*}
(X\tau)(F(t,q))
 &=d\tau_{F(t,q)}(F_*\partial_t)\\
 &=\partial_t(\tau(F(t,q)))=\partial_t t=1.
\end{align*}
\end{proof}

Any two admissible functions \(\tau\) differ by a smooth first
integral of \(X\), since
\[
X(\tau_1-\tau_2)=X\tau_1-X\tau_2=1-1=0.
\]
Along every complete orbit, each such function
takes every real value. In particular,
\(\operatorname{rank}(g-B)\) is zero on \(\tau^{-1}(0)\) and one
on its complement. Thus the rank in
Theorem~\ref{thm:complete-rank-one} cannot be required to be
identically one.

\begin{corollary}[Compactness obstruction]
\label{cor:compact-rank-one}
If \(M\) is compact and \(n\geq3\), no metric \(g\) and smooth
factor \(c\) satisfy the modified conformal equation together with
\(\operatorname{rank}(g-B)\leq1\). In dimension \(2\), the same
conclusion holds when the factor is required to be zero.
\end{corollary}

\begin{proof}
A smooth vector field on a compact manifold without boundary is
complete, by the continuation theorem for ordinary differential
equations. Theorem~\ref{thm:complete-rank-one} would therefore
give a smooth function \(\tau\) with \(X\tau=1\).
At a maximum of \(\tau\), compactness and the absence of boundary
give \(d\tau_p=0\), and hence
\[
1=(X\tau)(p)=d\tau_p(X_p)=0,
\]
a contradiction.
\end{proof}

\begin{example}[The dimension restriction]
\label{ex:dimension-two}
On \(\mathbb R^2\), with coordinates \((t,y)\), let
\[
 X=\partial_t,\qquad
 B=dt\otimes dt+dy\otimes dy,\qquad
 g=dt\otimes dt+e^{-2t/m}dy\otimes dy.
\]
The flow is
\[
\Phi_s(t,y)=(t+s,y),\qquad s\in\mathbb R,
\]
so \(X\) is complete. The matrices are
\[
[B]=\begin{pmatrix}1&0\\0&1\end{pmatrix},\quad
[g]=\begin{pmatrix}1&0\\0&e^{-2t/m}\end{pmatrix},\quad
[g^{-1}]=\begin{pmatrix}1&0\\0&e^{2t/m}\end{pmatrix}.
\]
Their inverse relation is verified by
\[
\begin{pmatrix}1&0\\0&e^{-2t/m}\end{pmatrix}
\begin{pmatrix}1&0\\0&e^{2t/m}\end{pmatrix}
=\begin{pmatrix}
1\cdot1+0\cdot0&1\cdot0+0\cdot e^{2t/m}\\
0\cdot1+e^{-2t/m}\cdot0&0\cdot0+e^{-2t/m}e^{2t/m}
\end{pmatrix}
=I_2.
\]
Both diagonal entries are positive, so the metrics are positive
definite. Since all coefficients of \(B\) are constant,
\[
\mathcal L_XB=\partial_t(dt\otimes dt+dy\otimes dy)=0.
\]
The difference
\[
 g-B=(e^{-2t/m}-1)dy\otimes dy
\]
has rank zero when \(t=0\) and rank one when \(t\ne0\).
Its contraction and derivative are
\begin{align*}
g(X,\cdot)
 &=dt(\partial_t)dt+e^{-2t/m}dy(\partial_t)dy
 =1\cdot dt+e^{-2t/m}0\cdot dy=dt,\\
\mathcal L_Xg
 &=\partial_t(e^{-2t/m})dy\otimes dy
 =-\frac2m e^{-2t/m}dy\otimes dy.
\end{align*}
Substitution into the modified conformal equation gives
\begin{align*}
\frac12\mathcal L_Xg-\frac1m
g(X,\cdot)\otimes g(X,\cdot)
 &=-\frac1m e^{-2t/m}dy\otimes dy-\frac1m dt\otimes dt\\
 &=-\frac1m(dt\otimes dt+e^{-2t/m}dy\otimes dy)
 =-\frac1m g.
\end{align*}
Thus the factor need not vanish in dimension \(2\), even for a
complete vector field. Since \(B(X,X)=1\), the global criterion
also fails in that dimension if nonzero factors are allowed.

The same formulas extend to \(\mathbb R^n\):
\[
 B=dt\otimes dt+\sum_{j=1}^{n-1}dy_j\otimes dy_j,\qquad
 g=dt\otimes dt+e^{-2t/m}\sum_{j=1}^{n-1}dy_j\otimes dy_j.
\]
The inverse is
\[
g^{-1}=\partial_t\otimes\partial_t+
e^{2t/m}\sum_{j=1}^{n-1}\partial_{y_j}\otimes\partial_{y_j}.
\]
In matrix form,
\[
[g][g^{-1}]
=\operatorname{diag}(1,e^{-2t/m},\ldots,e^{-2t/m})
 \operatorname{diag}(1,e^{2t/m},\ldots,e^{2t/m})
=\operatorname{diag}(1,1,\ldots,1).
\]
For \(X=\partial_t\), the computations are
\begin{align*}
g(X,\cdot)
 &=dt(\partial_t)dt+
 e^{-2t/m}\sum_{j=1}^{n-1}dy_j(\partial_t)dy_j=dt,\\
\frac12\mathcal L_Xg
 &=\frac12\partial_t(e^{-2t/m})
 \sum_{j=1}^{n-1}dy_j\otimes dy_j\\
 &=-\frac1m e^{-2t/m}\sum_{j=1}^{n-1}dy_j\otimes dy_j,\\
\frac12\mathcal L_Xg-\frac1m g(X,\cdot)\otimes g(X,\cdot)
 &=-\frac1m e^{-2t/m}\sum_{j=1}^{n-1}dy_j\otimes dy_j
 -\frac1m dt\otimes dt\\
 &=-\frac1m\left(dt\otimes dt+
 e^{-2t/m}\sum_{j=1}^{n-1}dy_j\otimes dy_j\right)
 =-\frac1m g.
\end{align*}
However,
\[
g-B=(e^{-2t/m}-1)\sum_{j=1}^{n-1}dy_j\otimes dy_j.
\]
For \(t\ne0\), all \(n-1\) spatial diagonal entries are nonzero
and the \(t\) entry is zero, so
\[
\operatorname{rank}(g-B)=n-1.
\]
For \(n\geq3\), the rank-one hypothesis is therefore not satisfied.
\end{example}

\section{Correspondences from an auxiliary vector field}
\label{sec:auxiliary-field}

In this section, \(X\) is nowhere zero and \(Z\) is a prescribed
smooth vector field satisfying
\begin{equation}
 [X,Z]=X.
 \label{eq:auxiliary-bracket}
\end{equation}
Completeness of \(X\) is not assumed.
We compare Killing metrics \(B\) with
\(B(X,X)=B(X,Z)=0\) and modified conformal metrics \(g\) with
factor zero and \(g(X,X)=g(X,Z)=0\).
Either class may be empty.
For a related shear construction in the Lie-algebraic setting,
see \cite[Lemma~4.6]{LDZ}.

\begin{theorem}[Correspondence under orthogonality conditions]
\label{thm:auxiliary-correspondence}
Suppose that a metric \(B\) satisfies
\[
 \mathcal L_XB=0,\qquad B(X,X)=B(X,Z)=0.
\]
Then the tensor
\begin{equation}
 g^{-1}=B^{-1}-\frac1m(X\otimes Z+Z\otimes X)
 \label{eq:auxiliary-inverse}
\end{equation}
is nondegenerate. Its inverse is
\begin{align}
 g={}&B+\frac1m\bigl(
 B(X,\cdot)\otimes B(Z,\cdot)
 +B(Z,\cdot)\otimes B(X,\cdot)\bigr)\notag\\
 &+\frac{B(Z,Z)}{m^2}
 B(X,\cdot)\otimes B(X,\cdot).
 \label{eq:auxiliary-covariant}
\end{align}
The metric \(g\) satisfies
\[
 \frac12\mathcal L_Xg
 =\frac1m g(X,\cdot)\otimes g(X,\cdot),\qquad
 g(X,X)=g(X,Z)=0.
\]
Conversely, every metric satisfying these three conditions gives
a Killing metric \(B\) with \(B(X,X)=B(X,Z)=0\) through
\[
 B^{-1}=g^{-1}+\frac1m(X\otimes Z+Z\otimes X).
\]
Explicitly,
\begin{align*}
 B={}&g-\frac1m\bigl(
 g(X,\cdot)\otimes g(Z,\cdot)
 +g(Z,\cdot)\otimes g(X,\cdot)\bigr)\\
 &+\frac{g(Z,Z)}{m^2}
 g(X,\cdot)\otimes g(X,\cdot).
\end{align*}
These maps are mutually inverse, preserve signature, and satisfy
\(g(X,\cdot)=B(X,\cdot)\).
At each point, \(\operatorname{rank}(g-B)\leq2\).
The rank is two if \(X,Z\) are linearly independent, one if
\(Z\) is a nonzero multiple of \(X\), and zero if \(Z=0\).
\end{theorem}

\begin{proof}
First consider the algebraic inversion. Regard
\[
 A=\frac1m\bigl(X\otimes B(Z,\cdot)+Z\otimes B(X,\cdot)\bigr)
\]
as an endomorphism of \(TM\). The elementary product rule
is obtained by evaluating on an arbitrary vector \(W\):
\begin{align*}
(U\otimes\mu)(V\otimes\nu)(W)
 &=(U\otimes\mu)(\nu(W)V)\\
 &=\mu(V)\nu(W)U
 =(\mu(V)\,U\otimes\nu)(W).
\end{align*}
Thus
\[
(U\otimes\mu)(V\otimes\nu)=\mu(V)\,U\otimes\nu.
\]
First,
\[
AX=\frac1m\bigl(B(Z,X)X+B(X,X)Z\bigr)
=\frac1m(0X+0Z)=0.
\]
Expanding all four terms in \(A^2\) gives
\begin{align*}
A^2=\frac1{m^2}\bigl[&
B(Z,X)X\otimes B(Z,\cdot)
+B(Z,Z)X\otimes B(X,\cdot)\\
&+B(X,X)Z\otimes B(Z,\cdot)
+B(X,Z)Z\otimes B(X,\cdot)\bigr]\\
=\frac1{m^2}\bigl[&
0X\otimes B(Z,\cdot)
+B(Z,Z)X\otimes B(X,\cdot)\\
&+0Z\otimes B(Z,\cdot)+0Z\otimes B(X,\cdot)\bigr]\\
= {}&\frac{B(Z,Z)}{m^2}X\otimes B(X,\cdot).
\end{align*}
Multiplication once more gives
\begin{align*}
A^3
 &=\frac{B(Z,Z)}{m^3}
 \bigl(B(Z,X)X\otimes B(X,\cdot)
 +B(X,X)Z\otimes B(X,\cdot)\bigr)\\
 &=\frac{B(Z,Z)}{m^3}
 \bigl(0X\otimes B(X,\cdot)+0Z\otimes B(X,\cdot)\bigr)=0.
\end{align*}
For every \(s\in\mathbb R\), the full product is
\begin{align*}
(I-sA)(I+sA+s^2A^2)
 &=I+sA+s^2A^2-sA-s^2A^2-s^3A^3\\
 &=I+sA-sA+s^2A^2-s^2A^2-0=I.
\end{align*}
The product in the other order is also \(I\), since both factors
are polynomials in \(A\). Hence
\[
(I-sA)^{-1}=I+sA+s^2A^2.
\]
As a map from \(T^*M\) to \(TM\), the tensor
\[
 B^{-1}-\frac{s}{m}(X\otimes Z+Z\otimes X)
 =(I-sA)B^{-1}
\]
is therefore invertible for every real \(s\). It is also symmetric.
The path from \(s=0\) to \(s=1\) consists of nondegenerate symmetric
forms, so its signature is constant.
At \(s=1\), inversion gives the full covariant expansion
\begin{align*}
g={}&B(I+A+A^2)\\
={}&B+\frac1m\bigl(B(X,\cdot)\otimes B(Z,\cdot)
 +B(Z,\cdot)\otimes B(X,\cdot)\bigr)\\
 &+\frac{B(Z,Z)}{m^2}B(X,\cdot)\otimes B(X,\cdot),
\end{align*}
which is~\eqref{eq:auxiliary-covariant}. Here the first line is
composition of maps, and the remaining lines are covariant tensors.
Thus \(B\) and \(g\) have the same signature.
For an arbitrary vector \(U\), its contractions are
\begin{align*}
g(X,U)
={}&B(X,U)+\frac1m
 \bigl(B(X,X)B(Z,U)+B(Z,X)B(X,U)\bigr)\\
 &+\frac{B(Z,Z)}{m^2}B(X,X)B(X,U)\\
={}&B(X,U)+\frac1m(0B(Z,U)+0B(X,U))
 +\frac{B(Z,Z)}{m^2}0B(X,U)\\
={}&B(X,U),\\
g(Z,U)
={}&B(Z,U)+\frac1m
 \bigl(B(X,Z)B(Z,U)+B(Z,Z)B(X,U)\bigr)\\
 &+\frac{B(Z,Z)}{m^2}B(X,Z)B(X,U)\\
={}&B(Z,U)+\frac{B(Z,Z)}m B(X,U),\\
g(Z,Z)
={}&B(Z,Z)+\frac1m
 \bigl(B(X,Z)B(Z,Z)+B(Z,Z)B(X,Z)\bigr)\\
 &+\frac{B(Z,Z)}{m^2}B(X,Z)^2
 =B(Z,Z).
\end{align*}
In particular,
\[
g(X,\cdot)=B(X,\cdot),\quad
g(X,X)=B(X,X)=0,\quad g(X,Z)=B(X,Z)=0.
\]

The differential equation follows directly from the inverse form.
Indeed, \([X,X]=0\) and~\eqref{eq:auxiliary-bracket} give
\begin{align*}
\mathcal L_X(X\otimes Z+Z\otimes X)
={}&[X,X]\otimes Z+X\otimes[X,Z]\\
 &+[X,Z]\otimes X+Z\otimes[X,X]\\
={}&0\otimes Z+X\otimes X+X\otimes X+Z\otimes0\\
={}&2X\otimes X.
\end{align*}
The inverse-metric identity gives
\[
\mathcal L_XB^{-1}
=-B^{-1}(\mathcal L_XB)B^{-1}
=-B^{-1}0B^{-1}=0.
\]
Consequently,
\begin{align*}
\mathcal L_Xg^{-1}
 &=\mathcal L_X\left(B^{-1}
 -\frac1m(X\otimes Z+Z\otimes X)\right)\\
 &=0-\frac1m(2X\otimes X)=-\frac2mX\otimes X.
\end{align*}
By Lemma~\ref{lem:inverse-equation}, \(g\) has modified conformal
factor zero.
One may also verify this in covariant form: the Killing condition
and the bracket identity give
\begin{align*}
\mathcal L_X(B(X,\cdot))
 &=(\mathcal L_XB)(X,\cdot)+B([X,X],\cdot)
 =0+B(0,\cdot)=0,\\
\mathcal L_X(B(Z,\cdot))
 &=(\mathcal L_XB)(Z,\cdot)+B([X,Z],\cdot)
 =0+B(X,\cdot)=B(X,\cdot),\\
X(B(Z,Z))
 &=(\mathcal L_XB)(Z,Z)+B([X,Z],Z)+B(Z,[X,Z])\\
 &=0+B(X,Z)+B(Z,X)=0+0+0=0.
\end{align*}
Applying these identities to~\eqref{eq:auxiliary-covariant} yields
\begin{align*}
\mathcal L_Xg
={}&0+\frac1m\bigl(
 0\otimes B(Z,\cdot)+B(X,\cdot)\otimes B(X,\cdot)\\
 &\hspace{24mm}+B(X,\cdot)\otimes B(X,\cdot)
 +B(Z,\cdot)\otimes0\bigr)\\
 &+\frac0{m^2}B(X,\cdot)\otimes B(X,\cdot)\\
 &+\frac{B(Z,Z)}{m^2}
 \bigl(0\otimes B(X,\cdot)+B(X,\cdot)\otimes0\bigr)\\
={}&\frac2m B(X,\cdot)\otimes B(X,\cdot)
 =\frac2m g(X,\cdot)\otimes g(X,\cdot).
\end{align*}

For the converse, start with a modified conformal metric \(g\) of
factor zero such that \(g(X,X)=g(X,Z)=0\).
The endomorphism
\[
 \widetilde A=\frac1m\bigl(
 X\otimes g(Z,\cdot)+Z\otimes g(X,\cdot)\bigr)
\]
has
\begin{align*}
\widetilde A^2=\frac1{m^2}\bigl[&
g(Z,X)X\otimes g(Z,\cdot)
+g(Z,Z)X\otimes g(X,\cdot)\\
&+g(X,X)Z\otimes g(Z,\cdot)
+g(X,Z)Z\otimes g(X,\cdot)\bigr]\\
=\frac1{m^2}\bigl[&
0X\otimes g(Z,\cdot)+g(Z,Z)X\otimes g(X,\cdot)\\
&+0Z\otimes g(Z,\cdot)+0Z\otimes g(X,\cdot)\bigr]\\
= {}&\frac{g(Z,Z)}{m^2}X\otimes g(X,\cdot),
\end{align*}
and
\begin{align*}
\widetilde A^3
 &=\frac{g(Z,Z)}{m^3}
 \bigl(g(Z,X)X\otimes g(X,\cdot)
 +g(X,X)Z\otimes g(X,\cdot)\bigr)\\
 &=\frac{g(Z,Z)}{m^3}
 \bigl(0X\otimes g(X,\cdot)+0Z\otimes g(X,\cdot)\bigr)=0.
\end{align*}
The inverse is checked by
\begin{align*}
&(I+s\widetilde A)(I-s\widetilde A+s^2\widetilde A^2)\\
&\qquad=I-s\widetilde A+s^2\widetilde A^2
 +s\widetilde A-s^2\widetilde A^2+s^3\widetilde A^3\\
&\qquad=I-s\widetilde A+s\widetilde A
 +s^2\widetilde A^2-s^2\widetilde A^2+0=I.
\end{align*}
The factors commute, so
\[
(I+s\widetilde A)^{-1}
=I-s\widetilde A+s^2\widetilde A^2.
\]
This proves that
\[
 B^{-1}=g^{-1}+\frac1m(X\otimes Z+Z\otimes X)
\]
is nondegenerate and has the same signature as \(g\), since the
path
\[
g^{-1}+\frac{s}{m}(X\otimes Z+Z\otimes X)
=(I+s\widetilde A)g^{-1}
\]
is symmetric and nondegenerate for every \(s\).
Inverting at \(s=1\) gives
\begin{align*}
B={}&g(I-\widetilde A+\widetilde A^2)\\
={}&g-\frac1m\bigl(g(X,\cdot)\otimes g(Z,\cdot)
 +g(Z,\cdot)\otimes g(X,\cdot)\bigr)\\
 &+\frac{g(Z,Z)}{m^2}g(X,\cdot)\otimes g(X,\cdot).
\end{align*}
Its contractions, for an arbitrary vector \(U\), are
\begin{align*}
B(X,U)
={}&g(X,U)-\frac1m
 \bigl(g(X,X)g(Z,U)+g(Z,X)g(X,U)\bigr)\\
 &+\frac{g(Z,Z)}{m^2}g(X,X)g(X,U)\\
={}&g(X,U)-\frac1m(0g(Z,U)+0g(X,U))
 +\frac{g(Z,Z)}{m^2}0g(X,U)=g(X,U),\\
B(Z,U)
={}&g(Z,U)-\frac1m
 \bigl(g(X,Z)g(Z,U)+g(Z,Z)g(X,U)\bigr)\\
 &+\frac{g(Z,Z)}{m^2}g(X,Z)g(X,U)\\
={}&g(Z,U)-\frac{g(Z,Z)}m g(X,U),\\
B(Z,Z)
={}&g(Z,Z)-\frac1m
 \bigl(g(X,Z)g(Z,Z)+g(Z,Z)g(X,Z)\bigr)\\
 &+\frac{g(Z,Z)}{m^2}g(X,Z)^2=g(Z,Z).
\end{align*}
In particular,
\[
B(X,\cdot)=g(X,\cdot),\quad
B(X,X)=g(X,X)=0,\quad B(X,Z)=g(X,Z)=0.
\]
The inverse-metric derivatives cancel:
\begin{align*}
\mathcal L_XB^{-1}
 &=\mathcal L_X\left(g^{-1}
 +\frac1m(X\otimes Z+Z\otimes X)\right)\\
 &=-\frac2mX\otimes X+\frac1m(2X\otimes X)=0.
\end{align*}
For an additional covariant verification, the required derivatives
are
\begin{align*}
\mathcal L_X(g(X,\cdot))
 &=(\mathcal L_Xg)(X,\cdot)+g([X,X],\cdot)\\
 &=\frac2m g(X,X)g(X,\cdot)+g(0,\cdot)=0,\\
\mathcal L_X(g(Z,\cdot))
 &=(\mathcal L_Xg)(Z,\cdot)+g([X,Z],\cdot)\\
 &=\frac2m g(X,Z)g(X,\cdot)+g(X,\cdot)=g(X,\cdot),\\
X(g(Z,Z))
 &=(\mathcal L_Xg)(Z,Z)+g([X,Z],Z)+g(Z,[X,Z])\\
 &=\frac2m g(X,Z)^2+g(X,Z)+g(Z,X)=0.
\end{align*}
Substitution into the covariant formula for \(B\) gives
\begin{align*}
\mathcal L_XB
={}&\frac2m g(X,\cdot)\otimes g(X,\cdot)\\
 &-\frac1m\bigl(
 0\otimes g(Z,\cdot)+g(X,\cdot)\otimes g(X,\cdot)\\
 &\hspace{24mm}+g(X,\cdot)\otimes g(X,\cdot)
 +g(Z,\cdot)\otimes0\bigr)\\
 &+\frac0{m^2}g(X,\cdot)\otimes g(X,\cdot)\\
 &+\frac{g(Z,Z)}{m^2}
 (0\otimes g(X,\cdot)+g(X,\cdot)\otimes0)\\
={}&\frac2m g(X,\cdot)\otimes g(X,\cdot)
 -\frac2m g(X,\cdot)\otimes g(X,\cdot)=0.
\end{align*}
Thus \(B\) is Killing.

The inverse translations cancel in both orders:
\begin{align*}
&\left[B^{-1}-\frac1m(X\otimes Z+Z\otimes X)\right]
 +\frac1m(X\otimes Z+Z\otimes X)=B^{-1},\\
&\left[g^{-1}+\frac1m(X\otimes Z+Z\otimes X)\right]
 -\frac1m(X\otimes Z+Z\otimes X)=g^{-1}.
\end{align*}
We also check the covariant substitutions.
For this calculation alone, put
\[
\beta=B(X,\cdot),\quad\gamma=B(Z,\cdot),\quad b=B(Z,Z).
\]
The established contractions read
\[
g(X,\cdot)=\beta,\qquad
g(Z,\cdot)=\gamma+\frac bm\beta,\qquad g(Z,Z)=b.
\]
Substituting these into the reverse formula gives
\begin{align*}
&g-\frac1m\bigl(g(X,\cdot)\otimes g(Z,\cdot)
 +g(Z,\cdot)\otimes g(X,\cdot)\bigr)
 +\frac{g(Z,Z)}{m^2}g(X,\cdot)\otimes g(X,\cdot)\\
&=B+\frac1m(\beta\otimes\gamma+\gamma\otimes\beta)
 +\frac b{m^2}\beta\otimes\beta\\
&\quad-\frac1m\left[
 \beta\otimes\left(\gamma+\frac bm\beta\right)
 +\left(\gamma+\frac bm\beta\right)\otimes\beta\right]
 +\frac b{m^2}\beta\otimes\beta\\
&=B+\frac1m(\beta\otimes\gamma+\gamma\otimes\beta)
 -\frac1m(\beta\otimes\gamma+\gamma\otimes\beta)\\
&\quad+\left(\frac b{m^2}-\frac b{m^2}
 -\frac b{m^2}+\frac b{m^2}\right)\beta\otimes\beta=B.
\end{align*}
In the other direction, put
\[
\widehat\beta=g(X,\cdot),\quad
\widehat\gamma=g(Z,\cdot),\quad\widehat b=g(Z,Z).
\]
Then
\[
B(X,\cdot)=\widehat\beta,\quad
B(Z,\cdot)=\widehat\gamma-\frac{\widehat b}m\widehat\beta,
\quad B(Z,Z)=\widehat b.
\]
Substitution into the forward formula gives
\begin{align*}
&B+\frac1m\bigl(B(X,\cdot)\otimes B(Z,\cdot)
 +B(Z,\cdot)\otimes B(X,\cdot)\bigr)
 +\frac{B(Z,Z)}{m^2}B(X,\cdot)\otimes B(X,\cdot)\\
&=g-\frac1m(\widehat\beta\otimes\widehat\gamma+
 \widehat\gamma\otimes\widehat\beta)
 +\frac{\widehat b}{m^2}\widehat\beta\otimes\widehat\beta\\
&\quad+\frac1m\left[
 \widehat\beta\otimes\left(\widehat\gamma-
 \frac{\widehat b}m\widehat\beta\right)
 +\left(\widehat\gamma-\frac{\widehat b}m\widehat\beta\right)
 \otimes\widehat\beta\right]
 +\frac{\widehat b}{m^2}\widehat\beta\otimes\widehat\beta\\
&=g-\frac1m(\widehat\beta\otimes\widehat\gamma+
 \widehat\gamma\otimes\widehat\beta)
 +\frac1m(\widehat\beta\otimes\widehat\gamma+
 \widehat\gamma\otimes\widehat\beta)\\
&\quad+\left(\frac{\widehat b}{m^2}
 -\frac{\widehat b}{m^2}-\frac{\widehat b}{m^2}
 +\frac{\widehat b}{m^2}\right)
 \widehat\beta\otimes\widehat\beta=g.
\end{align*}
Thus both the inverse and covariant formulas are mutually inverse.

For the rank statement, multiplication by invertible maps gives
\[
-g^{-1}(g-B)B^{-1}
=-g^{-1}gB^{-1}+g^{-1}BB^{-1}
=-B^{-1}+g^{-1}.
\]
Consequently,
\[
 \operatorname{rank}(g-B)
 =\operatorname{rank}(g^{-1}-B^{-1})
 =\operatorname{rank}\left(-\frac1m
 (X\otimes Z+Z\otimes X)\right).
\]
The image is contained in \(\operatorname{span}\{X,Z\}\).
If \(X,Z\) are independent, the matrix on their span, with respect
to these vectors and their dual covectors, is
\[
 \begin{pmatrix}0&-1/m\\-1/m&0\end{pmatrix},\qquad
 \det\begin{pmatrix}0&-1/m\\-1/m&0\end{pmatrix}
 =0\cdot0-\left(-\frac1m\right)\left(-\frac1m\right)
 =-\frac1{m^2}\ne0.
\]
Thus the rank is two.

If \(Z=\lambda X\) with \(\lambda\ne0\) at the point, then
\[
-\frac1m(X\otimes Z+Z\otimes X)
=-\frac1m(X\otimes\lambda X+\lambda X\otimes X)
=-\frac{2\lambda}mX\otimes X.
\]
This has rank one because \(X\ne0\).

If \(Z=0\), then
\[
-\frac1m(X\otimes0+0\otimes X)=0,
\]
so the rank is zero.
\end{proof}

\begin{example}[The bracket condition does not ensure nondegeneracy]
\label{ex:bracket-degeneracy}
On \(\mathbb R^2\), take
\[
 X=\partial_t,\qquad Z=t\partial_t+m\partial_y,\qquad
 B=dt\otimes dy+dy\otimes dt.
\]
The flow \(\Phi_s(t,y)=(t+s,y)\) is defined for all \(s\in\mathbb R\),
so \(X\) is complete. Constant coefficients give
\(\mathcal L_XB=0\). The bracket and independence are verified by
\[
[X,Z]=[\partial_t,t\partial_t+m\partial_y]
=(\partial_t t)\partial_t+(\partial_t m)\partial_y
=\partial_t=X,
\]
\[
\det\begin{pmatrix}1&t\\0&m\end{pmatrix}
=1\cdot m-0\cdot t=m\ne0.
\]
The contractions are
\begin{align*}
B(X,X)&=dt(\partial_t)dy(\partial_t)
 +dy(\partial_t)dt(\partial_t)=1\cdot0+0\cdot1=0,\\
B(X,Z)&=dt(\partial_t)dy(t\partial_t+m\partial_y)
 +dy(\partial_t)dt(t\partial_t+m\partial_y)\\
 &=1\cdot m+0\cdot t=m,\\
B(Z,Z)&=2\,dt(t\partial_t+m\partial_y)
 dy(t\partial_t+m\partial_y)=2tm.
\end{align*}
In the coordinate frame,
\[
 [B]=[B^{-1}]=
 \begin{pmatrix}0&1\\1&0\end{pmatrix},\qquad
 [X\otimes Z+Z\otimes X]=
 \begin{pmatrix}t&m\\0&0\end{pmatrix}
 +\begin{pmatrix}t&0\\m&0\end{pmatrix}
 =\begin{pmatrix}2t&m\\m&0\end{pmatrix}.
\]
For completeness,
\[
\det B=0\cdot0-1\cdot1=-1,\qquad
\begin{pmatrix}0&1\\1&0\end{pmatrix}
\begin{pmatrix}0&1\\1&0\end{pmatrix}
=\begin{pmatrix}0+1&0+0\\0+0&1+0\end{pmatrix}=I_2.
\]
Consequently,
\begin{align*}
\left[B^{-1}-\frac1m(X\otimes Z+Z\otimes X)\right]
 &=\begin{pmatrix}0&1\\1&0\end{pmatrix}
 -\frac1m\begin{pmatrix}2t&m\\m&0\end{pmatrix}\\
 &=\begin{pmatrix}-2t/m&1-m/m\\1-m/m&0\end{pmatrix}
 =\begin{pmatrix}-2t/m&0\\0&0\end{pmatrix},
\end{align*}
whose determinant is \((-2t/m)\cdot0-0\cdot0=0\).
It is therefore degenerate everywhere.
Thus the bracket identity, completeness, and \(B(X,X)=0\) do not
replace the condition \(B(X,Z)=0\).
\end{example}

For clarity, the exact nondegeneracy condition without either
orthogonality assumption is
\begin{align}
 &\det\left(B^{-1}-\frac1m(X\otimes Z+Z\otimes X)\right)\notag\\
 &\qquad=\det(B^{-1})
 \left[\left(1-\frac{B(X,Z)}m\right)^2
       -\frac{B(X,X)B(Z,Z)}{m^2}\right].
 \label{eq:auxiliary-determinant}
\end{align}
To verify the formula, choose a basis at the point in question,
put \(U=(X\ Z)\), and let \(V\) be the two-row matrix with rows
\(B(Z,\cdot)/m\) and \(B(X,\cdot)/m\).
Consider the block matrix
\[
 \begin{pmatrix}I_n&U\\V&I_2\end{pmatrix}
\]
and the elimination matrix
\(\begin{pmatrix}I_n&0\\-V&I_2\end{pmatrix}\), whose determinant
is one. Multiplication on the left gives
\[
\begin{pmatrix}I_n&0\\-V&I_2\end{pmatrix}
\begin{pmatrix}I_n&U\\V&I_2\end{pmatrix}
=\begin{pmatrix}I_n&U\\-V+V&-VU+I_2\end{pmatrix}
=\begin{pmatrix}I_n&U\\0&I_2-VU\end{pmatrix}.
\]
Multiplication on the right instead gives
\[
\begin{pmatrix}I_n&U\\V&I_2\end{pmatrix}
\begin{pmatrix}I_n&0\\-V&I_2\end{pmatrix}
=\begin{pmatrix}I_n-UV&U\\V-V&I_2\end{pmatrix}
=\begin{pmatrix}I_n-UV&U\\0&I_2\end{pmatrix}.
\]
Both products have the determinant of the original block matrix.
Taking their block triangular determinants yields
\[
\det(I_n)\det(I_2-VU)=\det(I_n-UV)\det(I_2),
\quad \det(I_n-UV)=\det(I_2-VU).
\]
The matrix product itself is
\[
VU=\frac1m\begin{pmatrix}
B(Z,X)&B(Z,Z)\\B(X,X)&B(X,Z)
\end{pmatrix}.
\]
Consequently, the determinant ratio is
\begin{align*}
&\frac{\det(B^{-1}-m^{-1}(X\otimes Z+Z\otimes X))}
 {\det(B^{-1})}\\
&\quad=\det\left[I-\frac1m
 \bigl(X\otimes B(Z,\cdot)+Z\otimes B(X,\cdot)\bigr)\right]\\
&\quad=\det(I_n-UV)=\det(I_2-VU)\\
&\quad=\det\begin{pmatrix}
1-B(Z,X)/m&-B(Z,Z)/m\\
-B(X,X)/m&1-B(X,Z)/m
\end{pmatrix}\\
&\quad=\left(1-\frac{B(Z,X)}m\right)
 \left(1-\frac{B(X,Z)}m\right)
 -\left(-\frac{B(Z,Z)}m\right)\left(-\frac{B(X,X)}m\right)\\
&\quad=\left(1-\frac{B(X,Z)}m\right)^2
 -\frac{B(X,X)B(Z,Z)}{m^2}\\
&\quad=1-\frac{2B(X,Z)}m
 +\frac{B(X,Z)^2}{m^2}-\frac{B(X,X)B(Z,Z)}{m^2}.
\end{align*}
This proves~\eqref{eq:auxiliary-determinant}.
When \(B\) is Killing and \([X,Z]=X\), the translated tensor
still satisfies the inverse equation with factor zero, because
\[
\mathcal L_X\left(B^{-1}-\frac1m(X\otimes Z+Z\otimes X)\right)
=0-\frac1m(2X\otimes X)=-\frac2mX\otimes X.
\]
It defines
a metric precisely where the scalar factor in
\eqref{eq:auxiliary-determinant} is nonzero.
In Example~\ref{ex:bracket-degeneracy}, this factor is
\[
\left(1-\frac mm\right)^2-\frac{0\cdot2tm}{m^2}
=(1-1)^2-0=0,
\]
in agreement with the matrix computation.

\begin{example}[The bracket condition does not ensure rank two]
\label{ex:bracket-rank-drop}
On \(\mathbb R^2\), let
\[
 X=\partial_t,\qquad Z=t\partial_t,\qquad
 B=dt\otimes dy+dy\otimes dt.
\]
The flow is again \(\Phi_s(t,y)=(t+s,y)\), so \(X\) is complete.
The bracket and contractions are
\begin{align*}
[X,Z]&=[\partial_t,t\partial_t]
 =(\partial_t t)\partial_t=\partial_t=X,\\
B(X,X)&=1\cdot0+0\cdot1=0,\\
B(X,Z)&=tB(X,X)=t\cdot0=0,\\
B(Z,Z)&=t^2B(X,X)=t^2\cdot0=0.
\end{align*}
Moreover,
\[
B(X,\cdot)=dy,\qquad B(Z,\cdot)=t\,dy.
\]
Substitution into Theorem~\ref{thm:auxiliary-correspondence} gives
\begin{align*}
g&=dt\otimes dy+dy\otimes dt+
\frac1m(dy\otimes t\,dy+t\,dy\otimes dy)
+\frac0{m^2}dy\otimes dy\\
&=dt\otimes dy+dy\otimes dt+\frac{2t}{m}dy\otimes dy,
\end{align*}
with
\[
 [g]=\begin{pmatrix}0&1\\1&2t/m\end{pmatrix},\qquad
 \det g=0\cdot(2t/m)-1\cdot1=-1,
\]
\[
 [g^{-1}]=\begin{pmatrix}-2t/m&1\\1&0\end{pmatrix}.
\]
The inverse is verified by
\[
\begin{pmatrix}0&1\\1&2t/m\end{pmatrix}
\begin{pmatrix}-2t/m&1\\1&0\end{pmatrix}
=\begin{pmatrix}
0(-2t/m)+1\cdot1&0\cdot1+1\cdot0\\
1(-2t/m)+(2t/m)1&1\cdot1+(2t/m)0
\end{pmatrix}
=I_2.
\]
Contraction and differentiation give
\begin{align*}
g(X,\cdot)
 &=dt(\partial_t)dy+dy(\partial_t)dt
 +\frac{2t}m dy(\partial_t)dy\\
 &=1\cdot dy+0\cdot dt+\frac{2t}m0\cdot dy=dy,\\
\mathcal L_Xg
 &=\partial_t(2t/m)dy\otimes dy
 =\frac2m dy\otimes dy.
\end{align*}
The modified conformal equation is therefore
\[
\frac12\mathcal L_Xg-\frac1m g(X,\cdot)\otimes g(X,\cdot)
=\frac1m dy\otimes dy-\frac1m dy\otimes dy=0.
\]
On the other hand,
\begin{align*}
g^{-1}-B^{-1}
 &=-\frac1m(\partial_t\otimes t\partial_t+
 t\partial_t\otimes\partial_t)\\
 &=-\frac{2t}m\partial_t\otimes\partial_t,
\end{align*}
whose coordinate matrix is
\[
\begin{pmatrix}-2t/m&0\\0&0\end{pmatrix}.
\]
It has rank one for \(t\ne0\) and rank zero for \(t=0\).
Completeness and the bracket identity alone therefore do not
imply rank two.
\end{example}

\section{Complete flows and constant-rank rigidity}\label{sec:constant-rank}

We now start with two given metrics: a Killing metric for $X$ and a
modified conformal metric with factor zero. Completeness converts the
evolution of their inverse-metric difference into a linear algebra problem
for a symmetric pencil defined for every real parameter.

\begin{lemma}[A symmetric pencil of constant rank]\label{lem:constant-rank-pencil}
Let $V$ be a finite-dimensional real vector space, let $0\ne x\in V$,
and let $S\in\operatorname{Sym}^2 V$ have rank $k$. Suppose that
\[
\operatorname{rank}(S+t\,x\otimes x)=k
\qquad\text{for every }t\in\mathbb R.
\]
Then $x\in\operatorname{im}S$. As an element of
$\operatorname{Sym}^2(\operatorname{im}S)$, the tensor $S$ is nondegenerate,
and its inverse satisfies
\[
\bigl(S|_{\operatorname{im}S}\bigr)^{-1}(x,x)=0.
\]
In particular, $k\ge2$. If $k=2$, there is a unique $z\in V$ such that
\[
S=-\frac1m(x\otimes z+z\otimes x).
\]
This vector $z$ is linearly independent of $x$.
\end{lemma}

Here $S|_{\operatorname{im}S}$ denotes the contravariant tensor supported
on $\operatorname{im}S$; it is not a restriction in the covariant sense.

\begin{proof}
Regard $S$ as a map $V^*\to V$. Its symmetry gives
\[
\lambda(S\mu)=\mu(S\lambda)
\qquad(\lambda,\mu\in V^*).
\]
In particular,
\begin{align*}
\lambda\in\ker S
&\Longleftrightarrow S\lambda=0\\
&\Longleftrightarrow \mu(S\lambda)=0
  \quad\text{for every }\mu\in V^*\\
&\Longleftrightarrow \lambda(S\mu)=0
  \quad\text{for every }\mu\in V^*\\
&\Longleftrightarrow \lambda|_{\operatorname{im}S}=0.
\end{align*}
Thus $\ker S$ is the annihilator of $\operatorname{im}S$.
Consequently, $S$ induces an isomorphism
$(\operatorname{im}S)^*\to\operatorname{im}S$, which proves the asserted
nondegeneracy.

If $x\notin\operatorname{im}S$, choose a basis of $V$ consisting first
of a basis of $\operatorname{im}S$, then $x$, and then any remaining
basis vectors. In this basis,
\[
[S]=\begin{pmatrix}A&0&0\\0&0&0\\0&0&0\end{pmatrix},
\qquad
[S+t\,x\otimes x]=
\begin{pmatrix}A&0&0\\0&t&0\\0&0&0\end{pmatrix},
\]
where $A$ is an invertible $k\times k$ matrix. For $t\ne0$,
\[
\operatorname{rank}(S+t\,x\otimes x)
=\operatorname{rank}A+\operatorname{rank}(t)=k+1\ne k.
\]
This contradicts the hypothesis. Hence $x\in\operatorname{im}S$.

Choose a basis of $\operatorname{im}S$ and let $A$ be the matrix of the
induced contravariant form. Multilinearity of the determinant gives
\begin{align*}
\det(A+t\,xx^{\mathsf T})
&=\det A+t\sum_{i,j=1}^k\operatorname{cof}_{ij}(A)x^ix^j\\
&=\det A+t\det A\sum_{i,j=1}^{k}(A^{-1})_{ji}x^ix^j\\
&=\det A\left(1+t\sum_{i,j=1}^{k}(A^{-1})_{ij}x^ix^j\right)\\
&=\det A\left(1+t\bigl(S|_{\operatorname{im}S}\bigr)^{-1}(x,x)\right).
\end{align*}
Indeed, every term involving two columns from $t\,xx^{\mathsf T}$
vanishes because these columns are proportional. Here
$\operatorname{cof}_{ij}(A)=(-1)^{i+j}\det A_{\widehat i\widehat j}$,
and the inverse formula gives
\[
\operatorname{cof}_{ij}(A)=\det A\,(A^{-1})_{ji}.
\]
For clarity, put $q=\bigl(S|_{\operatorname{im}S}\bigr)^{-1}(x,x)$
in the following scalar calculation. If $q\ne0$, the parameter
$t=-1/q$ gives
\[
\det(A-q^{-1}xx^{\mathsf T})
=\det A(1-q^{-1}q)=0.
\]
Since $S-q^{-1}x\otimes x$ is still supported on the $k$-dimensional
space $\operatorname{im}S$, it follows that
\[
\operatorname{rank}(S-q^{-1}x\otimes x)
=\operatorname{rank}(A-q^{-1}xx^{\mathsf T})<k,
\]
a contradiction. Thus $q=0$, which proves the nullity assertion.
The inclusions $0\ne x\in\operatorname{im}S$ give $k\ge1$. If $k=1$,
write $S=s\,e\otimes e$ and $x=\lambda e$ with $s\ne0$ and
$\lambda\ne0$. Then
\[
\bigl(S|_{\operatorname{im}S}\bigr)^{-1}(x,x)
=\frac{\lambda^2}{s}\ne0,
\]
again a contradiction. Therefore $k\ge2$.

Suppose that $k=2$. In a basis $(x,y)$ of $\operatorname{im}S$, write
\[
S=a\,x\otimes x+b(x\otimes y+y\otimes x)+d\,y\otimes y,
\qquad
[S]=\begin{pmatrix}a&b\\b&d\end{pmatrix}.
\]
Nondegeneracy gives $ad-b^2\ne0$. The inverse matrix is verified by
\[
\begin{pmatrix}a&b\\b&d\end{pmatrix}
\frac1{ad-b^2}\begin{pmatrix}d&-b\\-b&a\end{pmatrix}
=\frac1{ad-b^2}
\begin{pmatrix}ad-b^2&-ab+ab\\bd-bd&ad-b^2\end{pmatrix}
=I_2.
\]
Since $x$ has coordinates $(1,0)$, the nullity assertion gives
\[
0=\bigl(S|_{\operatorname{im}S}\bigr)^{-1}(x,x)
=\frac{d}{ad-b^2}.
\]
Consequently, $d=0$ and $b\ne0$. The vector
\[
z=-\frac{ma}{2}x-mb\,y
\]
satisfies
\begin{align*}
-\frac1m(x\otimes z+z\otimes x)
&=-\frac1m\left[
x\otimes\left(-\frac{ma}{2}x-mb\,y\right)
+\left(-\frac{ma}{2}x-mb\,y\right)\otimes x\right]\\
&=-\frac1m\left[-ma\,x\otimes x
-mb\,x\otimes y-mb\,y\otimes x\right]\\
&=a\,x\otimes x+b(x\otimes y+y\otimes x)=S.
\end{align*}
The coefficient $-mb$ of $y$ is nonzero, so $z$ is independent of $x$.

For uniqueness, suppose that $x\otimes v+v\otimes x=0$. Choose
$\theta\in V^*$ with $\theta(x)=1$. Contraction once and then twice
with $\theta$ gives
\begin{align*}
0&=(x\otimes v+v\otimes x)(\theta,\cdot)
  =\theta(x)v+\theta(v)x=v+\theta(v)x,\\
0&=\theta\bigl(v+\theta(v)x\bigr)
  =\theta(v)+\theta(v)\theta(x)=2\theta(v).
\end{align*}
Therefore $\theta(v)=0$ and $v=-\theta(v)x=0$. If $z_1,z_2$ both
give the asserted decomposition, subtracting the two formulas gives
\[
0=-\frac1m\bigl(x\otimes(z_1-z_2)+(z_1-z_2)\otimes x\bigr),
\]
so the preceding calculation with $v=z_1-z_2$ yields $z_1=z_2$.
\end{proof}

\begin{theorem}[The rank bound and the rank-two case]\label{thm:constant-rank-rigidity}
Let $X$ be nowhere zero and complete. Suppose that two given metrics
$B$ and $g$ satisfy
\[
\mathcal L_XB=0,
\qquad
\frac12\mathcal L_Xg
=\frac1m g(X,\cdot)\otimes g(X,\cdot).
\]
If $g^{-1}-B^{-1}$ has constant rank $k$ on $M$, then $k\ge2$.
Moreover, $k=2$ if and only if there is a unique smooth vector field
$Z$, pointwise linearly independent of $X$, such that
\begin{equation}\label{eq:rank-two-decomposition}
g^{-1}-B^{-1}=-\frac1m(X\otimes Z+Z\otimes X),
\qquad [X,Z]=X.
\end{equation}
In this case, for any smooth one-form $\theta$ with $\theta(X)=1$,
\begin{equation}\label{eq:extract-rank-two-field}
Z=-m(g^{-1}-B^{-1})(\theta,\cdot)
+\frac m2(g^{-1}-B^{-1})(\theta,\theta)X.
\end{equation}
The first contraction in this formula is a vector. In local coordinates,
its $i$th component is
$\sum_{j=1}^n(g^{ij}-B^{ij})\theta_j$.
\end{theorem}

\begin{proof}
The inverse-metric identity and the product rule give
\begin{align*}
0&=\mathcal L_X(g^{-1}g)
  =(\mathcal L_Xg^{-1})g+g^{-1}(\mathcal L_Xg),\\
\mathcal L_Xg^{-1}
&=-g^{-1}(\mathcal L_Xg)g^{-1}\\
&=-g^{-1}\left(\frac2m
  g(X,\cdot)\otimes g(X,\cdot)\right)g^{-1}.
\end{align*}
The contractions in the last expression are
\begin{align*}
(\mathcal L_Xg^{-1})^{ij}
&=-\frac2m\sum_{a,b=1}^n g^{ia}
 \left(\sum_{u=1}^n g_{au}X^u\right)
 \left(\sum_{v=1}^n g_{bv}X^v\right)g^{bj}\\
&=-\frac2m\sum_{u,v=1}^n
  \left(\sum_{a=1}^n g^{ia}g_{au}\right)
  \left(\sum_{b=1}^n g_{bv}g^{bj}\right)X^uX^v\\
&=-\frac2m\sum_{u,v=1}^n\delta^i_u\delta^j_vX^uX^v
 =-\frac2m X^iX^j.
\end{align*}
The Killing equation similarly gives
\[
\mathcal L_XB^{-1}=-B^{-1}(\mathcal L_XB)B^{-1}
=-B^{-1}(0)B^{-1}=0.
\]
Subtracting gives
\begin{align*}
\mathcal L_X(g^{-1}-B^{-1})
&=\mathcal L_Xg^{-1}-\mathcal L_XB^{-1}\\
&=-\frac2m X\otimes X-0=-\frac2m X\otimes X.
\end{align*}
Let $\Phi_t$ be the complete flow of $X$. A flow preserves its
generating field:
\[
(\Phi_t^*X)_p
=(d\Phi_t|_p)^{-1}X_{\Phi_t(p)}
=(d\Phi_t|_p)^{-1}(d\Phi_t|_p)X_p=X_p.
\]
Here the pullback of a contravariant tensor uses $(d\Phi_t)^{-1}$ in
each vector slot. Therefore
\begin{align*}
\frac{d}{dt}\bigl(\Phi_t^*(g^{-1}-B^{-1})\bigr)_p
&=\bigl(\Phi_t^*\mathcal L_X(g^{-1}-B^{-1})\bigr)_p\\
&=\left(\Phi_t^*\left(-\frac2m X\otimes X\right)\right)_p\\
&=-\frac2m(\Phi_t^*X)_p\otimes(\Phi_t^*X)_p\\
&=-\frac2m X_p\otimes X_p.
\end{align*}
Integration from $0$ to $t$ gives
\begin{align*}
\bigl(\Phi_t^*(g^{-1}-B^{-1})\bigr)_p
&=\bigl(\Phi_0^*(g^{-1}-B^{-1})\bigr)_p
+\int_0^t\left(-\frac2m X_p\otimes X_p\right)\,ds\\
&=(g^{-1}-B^{-1})_p-\frac{2t}{m}X_p\otimes X_p.
\end{align*}
Thus, for every $p\in M$ and $t\in\mathbb R$,
\begin{equation}\label{eq:inverse-difference-flow}
\bigl(\Phi_t^*(g^{-1}-B^{-1})\bigr)_p
=(g^{-1}-B^{-1})_p-\frac{2t}{m}X_p\otimes X_p.
\end{equation}
Pullback preserves rank. The constant-rank hypothesis therefore implies
\[
\operatorname{rank}\left((g^{-1}-B^{-1})_p
-\frac{2t}{m}X_p\otimes X_p\right)=k
\qquad(t\in\mathbb R).
\]
Completeness is essential here. For any $s\in\mathbb R$, substitute
$t=-ms/2$ to obtain
\begin{align*}
k&=\operatorname{rank}\left((g^{-1}-B^{-1})_p
-\frac{2(-ms/2)}m X_p\otimes X_p\right)\\
&=\operatorname{rank}\left((g^{-1}-B^{-1})_p+sX_p\otimes X_p\right).
\end{align*}
Lemma~\ref{lem:constant-rank-pencil}, applied with
$V=T_pM$, $x=X_p$, and $S=(g^{-1}-B^{-1})_p$, proves $k\ge2$.

If $k=2$, the same lemma gives a unique vector $Z_p$ at each point,
independent of $X_p$, satisfying the first equality in
\eqref{eq:rank-two-decomposition}. To prove smoothness, contract this
equality with a one-form $\theta$ such that $\theta(X)=1$:
\begin{align*}
(g^{-1}-B^{-1})(\theta,\cdot)
&=-\frac1m\bigl(\theta(X)Z+\theta(Z)X\bigr)\\
&=-\frac1m\bigl(Z+\theta(Z)X\bigr),\\
(g^{-1}-B^{-1})(\theta,\theta)
&=-\frac1m\bigl(\theta(X)\theta(Z)+\theta(Z)\theta(X)\bigr)\\
&=-\frac2m\theta(Z).
\end{align*}
Substituting both expressions into the right-hand side of
\eqref{eq:extract-rank-two-field} gives
\begin{align*}
&-m(g^{-1}-B^{-1})(\theta,\cdot)
+\frac m2(g^{-1}-B^{-1})(\theta,\theta)X\\
&\qquad=-m\left[-\frac1m\bigl(Z+\theta(Z)X\bigr)\right]
+\frac m2\left[-\frac2m\theta(Z)\right]X\\
&\qquad=Z+\theta(Z)X-\theta(Z)X=Z.
\end{align*}
Such a global smooth one-form exists: for any auxiliary positive
definite metric $H$, set
\[
\theta=\frac{H(X,\cdot)}{H(X,X)},\qquad
\theta(X)=\frac{H(X,X)}{H(X,X)}=1.
\]
Inserting this choice of $\theta$ gives the global formula
\begin{align*}
Z={}&-m(g^{-1}-B^{-1})
\left(\frac{H(X,\cdot)}{H(X,X)},\cdot\right)\\
&+\frac m2(g^{-1}-B^{-1})
\left(\frac{H(X,\cdot)}{H(X,X)},
\frac{H(X,\cdot)}{H(X,X)}\right)X.
\end{align*}
To make its smoothness explicit, write
\[
\theta_j=\frac{\sum_{a=1}^nH_{ja}X^a}{H(X,X)},
\qquad
H(X,X)=\sum_{a,b=1}^nH_{ab}X^aX^b>0.
\]
Then every component of $Z$ is
\begin{align*}
Z^i={}&-m\sum_{j=1}^n(g^{ij}-B^{ij})\theta_j
+\frac{mX^i}{2}\sum_{j,k=1}^n(g^{jk}-B^{jk})\theta_j\theta_k\\
={}&-\frac m{H(X,X)}
\sum_{j,a=1}^n(g^{ij}-B^{ij})H_{ja}X^a\\
&+\frac{mX^i}{2H(X,X)^2}
\sum_{j,k,a,b=1}^n(g^{jk}-B^{jk})H_{ja}H_{kb}X^aX^b.
\end{align*}
All numerators are smooth and the denominator is everywhere nonzero,
so $Z$ is smooth. Every admissible choice of $\theta$ gives the same
pointwise unique vector field.

Taking the Lie derivative of the decomposition gives
\begin{align*}
-\frac2m X\otimes X
&=\mathcal L_X(g^{-1}-B^{-1})\\
&=-\frac1m\mathcal L_X(X\otimes Z+Z\otimes X)\\
&=-\frac1m\bigl([X,X]\otimes Z+X\otimes[X,Z]\\
&\hspace{5.6em}+[X,Z]\otimes X+Z\otimes[X,X]\bigr)\\
&=-\frac1m\bigl(X\otimes[X,Z]+[X,Z]\otimes X\bigr),
\end{align*}
where $[X,X]=0$ was used in the last step.
Hence
\[
X\otimes([X,Z]-X)+([X,Z]-X)\otimes X=0.
\]
Two contractions with $\theta$ now give
\begin{align*}
0&=\theta(X)([X,Z]-X)+\theta([X,Z]-X)X\\
&=[X,Z]-X+\theta([X,Z]-X)X,\\
0&=\theta([X,Z]-X)+\theta([X,Z]-X)\theta(X)\\
&=2\theta([X,Z]-X).
\end{align*}
Consequently, $\theta([X,Z]-X)=0$, and the first equality reduces to
$[X,Z]-X=0$.

Conversely, if the decomposition in
\eqref{eq:rank-two-decomposition} holds with $X$ and $Z$ pointwise
independent, extend $(X,Z)$ to a basis at each point. The matrix of
the difference in that basis is
\[
\begin{pmatrix}
0&-1/m&0\\
-1/m&0&0\\
0&0&0
\end{pmatrix}.
\]
The leading block satisfies
\[
\det\begin{pmatrix}0&-1/m\\-1/m&0\end{pmatrix}
=0-\left(-\frac1m\right)\left(-\frac1m\right)
=-\frac1{m^2}\ne0.
\]
All remaining rows and columns vanish, so the rank is exactly two.
Finally, if $Z_1,Z_2$ yield the same decomposition, then
\[
X\otimes(Z_1-Z_2)+(Z_1-Z_2)\otimes X=0.
\]
Contraction gives
\begin{align*}
0&=(Z_1-Z_2)+\theta(Z_1-Z_2)X,\\
0&=\theta(Z_1-Z_2)+\theta(Z_1-Z_2)\theta(X)
 =2\theta(Z_1-Z_2).
\end{align*}
Thus $Z_1-Z_2=-\theta(Z_1-Z_2)X=0$, proving uniqueness.
\end{proof}

The converse concerns metrics $B$ and $g$ that are already given and
nondegenerate. The bracket relation $[X,Z]=X$ alone does not guarantee
that an inverse-metric translation remains nondegenerate. It also does
not imply $B(X,X)=B(X,Z)=0$. The correspondence in the preceding
section requires its stated orthogonality conditions.

\begin{example}[The constant-rank hypothesis]\label{ex:rank-jump}
On $\mathbb R^2$ with coordinates $(t,y)$, let
\[
X=\partial_t,\qquad
B=dt\otimes dy+dy\otimes dt,
\qquad
g=B+\frac{2t}{m}dy\otimes dy.
\]
The metric matrices and their inverses are
\[
[B]=[B^{-1}]=\begin{pmatrix}0&1\\1&0\end{pmatrix},
\quad
[g]=\begin{pmatrix}0&1\\1&2t/m\end{pmatrix},
\quad
[g^{-1}]=\begin{pmatrix}-2t/m&1\\1&0\end{pmatrix}.
\]
Their determinants and inverse products are
\begin{align*}
\det[B]&=0\cdot0-1\cdot1=-1,\\
\det[g]&=0\cdot(2t/m)-1\cdot1=-1,\\
[B][B^{-1}]
&=\begin{pmatrix}0&1\\1&0\end{pmatrix}
\begin{pmatrix}0&1\\1&0\end{pmatrix}
=\begin{pmatrix}1&0\\0&1\end{pmatrix},\\
[g][g^{-1}]
&=\begin{pmatrix}0&1\\1&2t/m\end{pmatrix}
\begin{pmatrix}-2t/m&1\\1&0\end{pmatrix}\\
&=\begin{pmatrix}1&0\\-2t/m+2t/m&1\end{pmatrix}=I_2.
\end{align*}
In particular, both metrics are nondegenerate. Since $X=\partial_t$,
the coordinate one-forms have zero Lie derivative. Hence
\begin{align*}
\mathcal L_XB
&=\partial_t(1)(dt\otimes dy+dy\otimes dt)=0,\\
\mathcal L_Xg
&=\mathcal L_XB+\partial_t(2t/m)dy\otimes dy
=\frac2m dy\otimes dy,\\
g(X,\cdot)
&=dt(\partial_t)dy+dy(\partial_t)dt
+\frac{2t}{m}dy(\partial_t)dy\\
&=dy+0+0=dy.
\end{align*}
Substitution into the modified conformal equation gives
\begin{align*}
\frac12\mathcal L_Xg
-\frac1m g(X,\cdot)\otimes g(X,\cdot)
&=\frac12\left(\frac2m dy\otimes dy\right)
-\frac1m dy\otimes dy\\
&=\frac1m dy\otimes dy-\frac1m dy\otimes dy=0.
\end{align*}
The flow is
\[
\Phi_s(t,y)=(t+s,y),\qquad
\Phi_0(t,y)=(t,y),\qquad
\frac d{ds}\Phi_s(t,y)=(1,0)=X_{\Phi_s(t,y)}.
\]
It is defined for every $s\in\mathbb R$, so $X$ is complete. However,
\begin{align*}
[g^{-1}-B^{-1}]
&=\begin{pmatrix}-2t/m&1\\1&0\end{pmatrix}
-\begin{pmatrix}0&1\\1&0\end{pmatrix}\\
&=\begin{pmatrix}-2t/m&0\\0&0\end{pmatrix},\\
g^{-1}-B^{-1}&=-\frac{2t}{m}\partial_t\otimes\partial_t.
\end{align*}
Therefore
\[
\operatorname{rank}(g^{-1}-B^{-1})_{(t,y)}
=\begin{cases}0,&t=0,\\1,&t\ne0.\end{cases}
\]
This shows why the
constant-rank hypothesis cannot be omitted from the rank bound.
\end{example}

\begin{example}[The completeness hypothesis]\label{ex:incomplete-rank-one}
On $M=(0,m/2)\times\mathbb R$, take
\[
X=\partial_t,\qquad
B=dt\otimes dt+dy\otimes dy,
\qquad
g=\frac{dt\otimes dt}{1-2t/m}+dy\otimes dy.
\]
The matrices are
\[
[B]=[B^{-1}]=I_2,\qquad
[g]=\begin{pmatrix}(1-2t/m)^{-1}&0\\0&1\end{pmatrix}.
\]
Since $0<t<m/2$, their diagonal entries and determinants satisfy
\[
1-\frac{2t}{m}>0,\qquad
\det[B]=1,\qquad
\det[g]=\frac1{1-2t/m}>0.
\]
Both metrics are therefore positive definite. The inverse is verified by
\begin{align*}
[B][B^{-1}]&=I_2I_2=I_2,\\
[g]\begin{pmatrix}1-2t/m&0\\0&1\end{pmatrix}
&=\begin{pmatrix}(1-2t/m)^{-1}&0\\0&1\end{pmatrix}
\begin{pmatrix}1-2t/m&0\\0&1\end{pmatrix}\\
&=\begin{pmatrix}(1-2t/m)^{-1}(1-2t/m)&0\\0&1\end{pmatrix}
=I_2.
\end{align*}
The derivative of the nonconstant coefficient is
\[
\frac d{dt}(1-2t/m)^{-1}
=-(1-2t/m)^{-2}\left(-\frac2m\right)
=\frac2{m(1-2t/m)^2}.
\]
Consequently,
\begin{align*}
\mathcal L_XB
&=\partial_t(1)dt\otimes dt+\partial_t(1)dy\otimes dy=0,\\
\mathcal L_Xg
&=\frac2{m(1-2t/m)^2}dt\otimes dt
+\partial_t(1)dy\otimes dy\\
&=\frac2{m(1-2t/m)^2}dt\otimes dt,\\
g(X,\cdot)
&=\frac{dt(\partial_t)}{1-2t/m}dt+dy(\partial_t)dy
=\frac{dt}{1-2t/m}.
\end{align*}
The quadratic term and the full equation are
\begin{align*}
\frac1m g(X,\cdot)\otimes g(X,\cdot)
&=\frac1m\left(\frac{dt}{1-2t/m}\right)
\otimes\left(\frac{dt}{1-2t/m}\right)\\
&=\frac{dt\otimes dt}{m(1-2t/m)^2},\\
\frac12\mathcal L_Xg
-\frac1m g(X,\cdot)\otimes g(X,\cdot)
&=\frac{dt\otimes dt}{m(1-2t/m)^2}
-\frac{dt\otimes dt}{m(1-2t/m)^2}=0.
\end{align*}
Subtracting the inverse matrices yields
\begin{align*}
[g^{-1}-B^{-1}]
&=\begin{pmatrix}1-2t/m&0\\0&1\end{pmatrix}
-\begin{pmatrix}1&0\\0&1\end{pmatrix}\\
&=\begin{pmatrix}-2t/m&0\\0&0\end{pmatrix}.
\end{align*}
Thus
\[
t>0\quad\Longrightarrow\quad -\frac{2t}{m}\ne0
\quad\Longrightarrow\quad\operatorname{rank}(g^{-1}-B^{-1})=1.
\]
The integral curve from $(t,y)$ satisfies
\[
\gamma(s)=(t+s,y),\qquad \dot\gamma(s)=(1,0)=X_{\gamma(s)}.
\]
Its maximal parameter interval in $M$ is determined by
\[
\gamma(s)\in M
\quad\Longleftrightarrow\quad 0<t+s<m/2
\quad\Longleftrightarrow\quad -t<s<m/2-t.
\]
This interval is not $\mathbb R$, so $X$ is not complete.
\end{example}

\begin{example}[Rank two need not preserve the dual one-form]\label{ex:rank-two-dual}
On $\mathbb R^2$, let
\[
X=\partial_t,\qquad Z=t\partial_t+m\partial_y,
\]
and set
\[
B=\frac12(dt\otimes dy+dy\otimes dt),\qquad
g=dt\otimes dy+dy\otimes dt+\frac{2t}{m}dy\otimes dy.
\]
Their matrices and inverse matrices are
\[
[B]=\begin{pmatrix}0&1/2\\1/2&0\end{pmatrix},\quad
[g]=\begin{pmatrix}0&1\\1&2t/m\end{pmatrix},
\]
\[
[B^{-1}]=\begin{pmatrix}0&2\\2&0\end{pmatrix},\quad
[g^{-1}]=\begin{pmatrix}-2t/m&1\\1&0\end{pmatrix}.
\]
The determinant calculations and inverse products are
\begin{align*}
\det[B]&=0-\frac12\cdot\frac12=-\frac14,\\
\det[g]&=0\cdot\frac{2t}{m}-1\cdot1=-1,\\
[B][B^{-1}]
&=\begin{pmatrix}0&1/2\\1/2&0\end{pmatrix}
\begin{pmatrix}0&2\\2&0\end{pmatrix}
=\begin{pmatrix}(1/2)2&0\\0&(1/2)2\end{pmatrix}=I_2,\\
[g][g^{-1}]
&=\begin{pmatrix}0&1\\1&2t/m\end{pmatrix}
\begin{pmatrix}-2t/m&1\\1&0\end{pmatrix}\\
&=\begin{pmatrix}1&0\\-2t/m+2t/m&1\end{pmatrix}=I_2.
\end{align*}
The negative determinants show that both metrics are nondegenerate
with signature $(1,1)$. Their inverse difference is
\begin{align*}
[g^{-1}-B^{-1}]
&=\begin{pmatrix}-2t/m&1\\1&0\end{pmatrix}
-\begin{pmatrix}0&2\\2&0\end{pmatrix}\\
&=\begin{pmatrix}-2t/m&-1\\-1&0\end{pmatrix},\\
\det[g^{-1}-B^{-1}]
&=\left(-\frac{2t}{m}\right)0-(-1)(-1)=-1\ne0.
\end{align*}
Therefore $\operatorname{rank}(g^{-1}-B^{-1})=2$ everywhere.
Substituting the prescribed vector fields into the decomposition gives
\begin{align*}
-\frac1m(X\otimes Z+Z\otimes X)
&=-\frac1m\bigl[
\partial_t\otimes(t\partial_t+m\partial_y)
+(t\partial_t+m\partial_y)\otimes\partial_t\bigr]\\
&=-\frac1m\bigl[2t\partial_t\otimes\partial_t
+m\partial_t\otimes\partial_y+m\partial_y\otimes\partial_t\bigr]\\
&=-\frac{2t}{m}\partial_t\otimes\partial_t
-\partial_t\otimes\partial_y-\partial_y\otimes\partial_t\\
&=g^{-1}-B^{-1}.
\end{align*}
The bracket is
\begin{align*}
[X,Z]&=[\partial_t,t\partial_t+m\partial_y]\\
&=(\partial_t t)\partial_t+(\partial_t m)\partial_y
+t[\partial_t,\partial_t]+m[\partial_t,\partial_y]\\
&=\partial_t+0+0+0=X.
\end{align*}
For every real $s$, the flow satisfies
\[
\Phi_s(t,y)=(t+s,y),\qquad
\frac d{ds}\Phi_s(t,y)=(1,0)=X_{\Phi_s(t,y)}.
\]
It is defined on all of $\mathbb R$, so $X$ is complete.
The metric derivatives are
\begin{align*}
\mathcal L_XB
&=\partial_t(1/2)(dt\otimes dy+dy\otimes dt)=0,\\
\mathcal L_Xg
&=\partial_t(1)(dt\otimes dy+dy\otimes dt)
+\partial_t(2t/m)dy\otimes dy\\
&=\frac2m dy\otimes dy.
\end{align*}
The dual one-forms are
\begin{align*}
B(X,\cdot)
&=\frac12\bigl(dt(\partial_t)dy+dy(\partial_t)dt\bigr)
=\frac12(1\,dy+0\,dt)=\frac12dy,\\
g(X,\cdot)
&=dt(\partial_t)dy+dy(\partial_t)dt
+\frac{2t}{m}dy(\partial_t)dy\\
&=dy+0+0=dy.
\end{align*}
Inserting these formulas gives
\begin{align*}
\frac12\mathcal L_Xg
-\frac1m g(X,\cdot)\otimes g(X,\cdot)
&=\frac12\left(\frac2m dy\otimes dy\right)
-\frac1m dy\otimes dy\\
&=\frac1m dy\otimes dy-\frac1m dy\otimes dy=0.
\end{align*}
Thus $g$ has modified conformal factor zero, whereas
$B(X,\cdot)=\tfrac12dy\ne dy=g(X,\cdot)$. Finally,
\begin{align*}
B(X,Z)
&=\frac12dy(t\partial_t+m\partial_y)
=\frac12\bigl(t\,dy(\partial_t)+m\,dy(\partial_y)\bigr)
=\frac m2,\\
g(X,Z)
&=dy(t\partial_t+m\partial_y)
=t\,dy(\partial_t)+m\,dy(\partial_y)=m.
\end{align*}
The orthogonality conditions from the preceding section therefore fail.
\end{example}

\begin{example}[The zero-factor hypothesis]\label{ex:nonzero-factor-rank-one}
On $\mathbb R^2$, take
\[
X=\partial_t,\qquad
B=dt\otimes dt+dy\otimes dy,
\qquad
g=dt\otimes dt-e^{-2t/m}dy\otimes dy.
\]
The metric matrices and inverse matrices are
\[
[B]=[B^{-1}]=I_2,\qquad
[g]=\begin{pmatrix}1&0\\0&-e^{-2t/m}\end{pmatrix},
\qquad
[g^{-1}]=\begin{pmatrix}1&0\\0&-e^{2t/m}\end{pmatrix}.
\]
Their determinants are
\[
\det[B]=1,\qquad
\det[g]=1\cdot(-e^{-2t/m})-0=-e^{-2t/m}\ne0.
\]
The displayed inverses are verified by
\begin{align*}
[B][B^{-1}]&=I_2I_2=I_2,\\
[g][g^{-1}]
&=\begin{pmatrix}1&0\\0&-e^{-2t/m}\end{pmatrix}
\begin{pmatrix}1&0\\0&-e^{2t/m}\end{pmatrix}\\
&=\begin{pmatrix}1&0\\0&(-e^{-2t/m})(-e^{2t/m})\end{pmatrix}
=\begin{pmatrix}1&0\\0&e^{-2t/m+2t/m}\end{pmatrix}=I_2.
\end{align*}
The flow and its derivative are
\[
\Phi_s(t,y)=(t+s,y),\qquad
\frac d{ds}\Phi_s(t,y)=(1,0)=X_{\Phi_s(t,y)}
\quad(s\in\mathbb R),
\]
so $X$ is complete. The metric derivatives and the dual one-form are
\begin{align*}
\mathcal L_XB
&=\partial_t(1)dt\otimes dt+\partial_t(1)dy\otimes dy=0,\\
\mathcal L_Xg
&=\partial_t(1)dt\otimes dt
+\partial_t(-e^{-2t/m})dy\otimes dy\\
&=-\left(-\frac2m e^{-2t/m}\right)dy\otimes dy
=\frac2m e^{-2t/m}dy\otimes dy,\\
g(X,\cdot)
&=dt(\partial_t)dt-e^{-2t/m}dy(\partial_t)dy
=dt-0=dt.
\end{align*}
Therefore
\begin{align*}
\frac12\mathcal L_Xg
-\frac1m g(X,\cdot)\otimes g(X,\cdot)
&=\frac1m e^{-2t/m}dy\otimes dy-\frac1m dt\otimes dt\\
&=-\frac1m\bigl(dt\otimes dt-e^{-2t/m}dy\otimes dy\bigr)
=-\frac1m g.
\end{align*}
The modified conformal factor is the nonzero constant $-1/m$.
Nevertheless,
\begin{align*}
[g^{-1}-B^{-1}]
&=\begin{pmatrix}1&0\\0&-e^{2t/m}\end{pmatrix}
-\begin{pmatrix}1&0\\0&1\end{pmatrix}\\
&=\begin{pmatrix}0&0\\0&-(e^{2t/m}+1)\end{pmatrix},\\
g^{-1}-B^{-1}&=-(e^{2t/m}+1)\partial_y\otimes\partial_y.
\end{align*}
Since $e^{2t/m}+1>0$, the only nonzero diagonal entry is everywhere
nonzero, and hence
\[
\operatorname{rank}(g^{-1}-B^{-1})=1
\qquad\text{on all of }\mathbb R^2.
\]
Thus the zero-factor hypothesis in
Theorem~\ref{thm:constant-rank-rigidity} cannot be replaced by an
arbitrary constant factor. This example uses metrics of different
signatures, which are permitted in that theorem.
\end{example}


\end{document}